\documentstyle[11pt]{article}
\newtheorem{thm}{Theorem}[section]
\newtheorem{lem}[thm]{Lemma}
\newtheorem{cor}[thm]{Corollary}
\newtheorem{pro}[thm]{Proposition}
\newtheorem{ex}[thm]{Example}

\newtheorem{defi}[thm]{Definition}

\newcommand{\gm }{\Gamma }
\newcommand{\lon }{\longrightarrow }
\newcommand{\be }{\begin{eqnarray*}}
\newcommand{\ee }{\end{eqnarray*}}

\newcommand{\poidd }[2]{#1\gpd #2}
\newcommand{\poiddd }[3]{ (#1\gpd #2, \alpha_{#3}, \beta_{#3})}

\input amssym.def
\input amssym

\newcommand{\pf}{\noindent{\bf Proof.}\ }
\newcommand{\qed}{\begin{flushright} {\bf Q.E.D.}\ \ \ \ \ \ \end{flushright}}

\newcommand{\reals}{{\Bbb  R}}

\newcommand{\frakg}{{\frak g}}
\newcommand{\frakh}{{\frak h}}

\newcommand{\half}{\frac{1}{2}}
\newcommand{\third}{\frac{1}{3}}

\newcommand{\D}{{\cal D}}

\newcommand{\cald}{{\cal D}}

\newcommand{\calf}{{\cal F}}

\newcommand{\calr}{{\cal R}}
\newcommand{\calx}{{\cal X}}
\newcommand{\caly}{{\cal Y}}

\newcommand{\smalcirc}{\mbox{\tiny{$\circ $}}}

\def\description label#1{\hfil\bf[#1]\hfil}
\newcommand{\kkk}{\phi_{k*}}
\newcommand{\kk}{\phi_{k}}
\newcommand{\F}{\phi }
\newcommand{\G}{\psi}
\newcommand{\closed}{ regular }
\newcommand{\Closed}{ Regular }
\newcommand{\regular}{ reducible }

\def\sdp{\mathbin{\hbox{$\mapstochar\kern-.3333em\times$}}}
\def\pds{\mathbin{\hbox{$\times\kern-.55em\mapstochar\,$}}}

\newcommand{\wed}{\mathbin{\lower1.5pt\hbox{$\scriptstyle{\wedge}$}}}

\let\Tilde=\widetilde

\def\chigh{{\raise1.5pt\hbox{$\chi$}}}
\let\phi=\varphi

\def\til0{\Tilde{0}}

\def\dminus{\raise2pt\hbox{\vrule height1pt width 2ex}\hskip3pt}

\def\pback#1{\mathbin{{{\lower1.2ex\hbox{$\times$}}\atop #1}}}

\def\vlra{\hbox{$\,-\!\!\!-\!\!\!-\!\!\!-\!\!\!-\!\!\!
-\!\!\!-\!\!\!-\!\!\!-\!\!\!-\!\!\!\longrightarrow\,$}}

\def\gpd{\,\lower1pt\hbox{$\longrightarrow$}\hskip-.24in\raise2pt
             \hbox{$\longrightarrow$}\,}

\def\lgpd{\,\lower1pt\hbox{$\vlra$}\hskip-1.02in\raise2pt\hbox{$\vlra$}\,}

\def\llgpd{\,\lower1pt\hbox{$\vvlra$}\hskip-1.3in\raise2pt\hbox{$\vvlra$}\,}

\begin{document}

\title{{\bf Dirac structures and Poisson homogeneous spaces}
}

\author{ ZHANG-JU LIU\thanks{Research supported by NSF of China.}\\
        Department  of Mathematics \\
        Peking University \\
        Beijing, 100871, China \\
        {\sf email: liuzj@sxx0.math.pku.edu.cn}\\
        ALAN WEINSTEIN \thanks{Research partially supported by NSF
        grants DMS93-09653 and DMS96-25122.}\\
        Department of  Mathematics\\
        University of California\\
        Berkeley, CA 94720, USA\\
        {\sf email: alanw@math.berkeley.edu}\\
        and \\
        PING XU \thanks {Research partially supported by NSF
        Grant DMS95-04913.}\\
Department of Mathematics\\
The  Pennsylvania State University\\
 University Park, PA 16802, USA\\
        {\sf email: ping@math.psu.edu }}

 \date{September 1996}

\maketitle
\begin{abstract}
Poisson homogeneous spaces for Poisson groupoids are classfied in
terms of Dirac structures for the corresponding Lie bialgebroids.
Applications include Drinfel'd's classification in the case of
Poisson groups and a description of leaf spaces of foliations as
homogeneous spaces of pair groupoids.

\end{abstract}
\section{Introduction}
Dirac structures on manifolds include closed 2-forms, Poisson
structures, and foliations.  They extend the flexibility of
computations with such objects by permitting the passage to both
submanifolds and quotients.  The combination of these two operations
is central to the theory of reduction in Poisson geometry.  Dirac
structures were introduced by Courant and Weinstein
\cite{co-we:beyond} and thoroughly investigated by Courant in
\cite{Courant:1990}.  Dorfman \cite{Dorfman} used Dirac structures in
the context of the formal calculus of variations  for the study of
completely integrable systems of partial differential equations.

Under a regularity assumption which is always satisfied on an open
dense subset, a Dirac structure on a manifold $P$ is locally the same
thing as a Poisson structure on the leaf space of a foliation of
$P$.\footnote{Dually, under a slightly different regularity
assumption, it is the same as a smooth family of closed 2-forms on the
leaves of a foliation of $P$ (generally different from the foliation
in the first description).}   

An essential object in the theory of Dirac structures 
is a natural antisymmetric bracket operation (see Equation
(\ref{eq:courant})) on the sections of $TP\oplus T^*P$ introduced by
Courant.  Although this Courant bracket does not satisfy the Jacobi
identity, it does satisfy that identity on $\Gamma(E)$ when $E$ is a
subbundle of $TP\oplus T^*P$ which is maximal isotropic for the
symmetric form $(X_{1}+\xi_{1} , X_{2}+\xi_{2})_{+}=\half (\langle
\xi_{1},  X_{2}  \rangle
  + \langle \xi_{2} ,  X_{1}\rangle )$ and whose sections are closed
under the bracket.  (For instance, we recover the usual bracket of
vector fields on $TP\oplus {0}$ and the zero bracket on ${0}\oplus
T^*P$.)  

The theory of Dirac structures finds an echo in Drinfel'd's theory of
Lie bialgebras and Poisson homogeneous spaces
\cite{dr:hamiltonian}\cite{dr:poisson}.  A Lie bialgebra $(\frakg ,
\frakg^{*})$ can be
thought of as a pair of Lie algebra structures on a vector space
$\frakg$ and its dual having a common extension (which turns out to be
unique) to a Lie algebra structure on $\frakg\oplus\frakg^*$ for which
the symmetric form $(~,~)_+$ is ad-invariant.  The Lie algebra
$\frakg\oplus\frakg^* $ is called the {\bf double} of the Lie bialgebra 
 $(\frakg ,\frakg^{*})$.  

The main result of \cite{dr:poisson} is that maximal isotropic
subalgebras of the double correspond (modulo some details concerning
closedness and connectedness of subgroups) to Poisson homogeneous
$G$-spaces, where $G$ is the Poisson Lie group whose linearization is
the given Lie bialgebra.

The similarities between the Courant bracket and the bracket on the double
of a Lie bialgebra were explained in our recent paper \cite{LWX},
where both were exhibited as special cases in a theory of doubles of
Lie bialgebr{\em oid}s.  Since the bracket on sections of a Lie algebroid
$A$ over $P$ satisfies the Jacobi identity, while the Courant bracket does not,
it is clear that one must look beyond Lie algebroids to find these
doubles.  Hence we introduced in \cite{LWX} a notion of {\bf Courant
algebroid}, in which the Lie algebroid axioms for a bracket on
$\Gamma(A)$ and a bundle map $a:A\to P$ are satisfied only modulo
certain ``coboundary anomalies,'' explicitly described in terms of a
nondegenerate bilinear form on $E$ which is part of the Courant
algebroid structure.  When $P$ is a point, the anomalies vanish, and a
Courant algebroid is just a Lie algebra with a nondegenerate
ad-invariant symmetric bilinear form.  

Although more explicit descriptions are available (see Section
\ref{sec:dirac}), we can define a Lie bialgebroid as a pair of Lie
algebroid structures on a vector bundle $A$ and its dual having a
common extension (which turns out to be unique) to a Courant algebroid
structure on $A\oplus A^*$ with the symmetric form $(~,~)_+$.  For the
original Courant bracket, $A=TP$ with the usual bracket of vector
fields and $A^*=T^*P$ with the zero bracket on 1-fomrs.  When $P$ is a
point, we recover the Lie bialgebras.  Lie bialgebroids arise as the
linearizations of (possible local) Poisson groupoids $G$, in
which the bracket on $A$ determines $G$,
while the bracket on $A^*$ determines a compatible Poisson structure.

If $(A,A^*)$ is a Lie bialgebroid over $P$, we can now define an
$(A,A^*)$ {\bf Dirac structure on $P$} to be a maximal isotropic subbundle
of the Courant algebroid $L\subset A\oplus A^*$ which is closed under
the bracket.  Since the Lie algebroid anomalies are defined in terms
of $(~,~)_+$, they vanish on $L$, which is therefore an ordinary Lie
algebroid.  

The main result of this paper, already announced in \cite{LWX}, is
that there is a 1-1 correspondence between $(A,A^*)$ Dirac structures
on $P$ satisfying a certain regularity condition and Poisson
homogeneous spaces of the form $G/H$, where $G$ is a Poisson groupoid
whose tangent Lie bialgebroid is $(A,A^*)$, and $H$ is a subgroupoid
of $G$ which is closed and wide, i.e. containing all the identity
elements.  (We will assume throughout this paper that all our
groupoids 
are $\alpha$-connected; i.e. the fibres of the source map
are connected).  Drinfel'd's theorem is the special case of this for
$P$ a point, while for ordinary (i.e. $(TP,T^*P)$) Dirac structures, we
recover their description as Poisson structures on quotient manifolds
of $P$.

On the way to our main result, we develop several topics of
independent interest.  First of all, we extend to $(A,A^*)$ Dirac
structures the original
application of Dirac structures to Poisson reduction.  A technical
complication here is that we must deal with quotient spaces for which
the projection is {\em not} a Poisson mapping.  This renders the
streamlined methods of ``coisotropic calculus'' \cite{we:coisotropic}
inapplicable, and we need to do a number of computations by hand.
Eventually, it will be useful to develop a modified coisotropic calculus
to hand Dirac reductions directly.
Second, we study pullbacks of Dirac structures under morphisms of Lie
bialgebroids.  Finally, we discuss the general notion of homogeneous
spaces for groupoids.

Here is an outline of the paper.  Section 2 is a review of basic
definitions and properties of Lie bialgebroids, Courant algebroids,
and Dirac structures (we will often omit the prefix ``$(A,A^*)$", which
is still implied).  In Section 3, we establish a correspondence between
Dirac structures and Poisson structures on quotient manifolds.  Using
this correspondence, we characterize in Section 4 the Dirac structures
which are invariant under Poisson actions of groups, and we prove
Drinfel'd's theorem by extending a $(\frakg,\frakg*)$ Dirac structure
 to a left-invariant $(TG,T^*G)$ Dirac structure on $G$.
(Here, the bialgebroid has the nontrivial bracket on 1-forms coming
from the Poisson structure on $G$.)  Section 5 contains a theorem
about pullbacks of Dirac structures, which is used in Section 6 to
extend $(A,A^* )$ Dirac structures on $P$ to ``left-invariant''
$(TG,T^*G)$ Dirac structures on the Poisson groupoid $G$.  These
invariant structures are then related to Poisson structures on
quotients of $P$.  Section 7 establishes a characterization of Poisson
actions of Poisson groupoids.  Finally, in Section 8, we define 
Poisson homogeneous spaces and then use the results of
Section 6 to prove our main theorem.

{\bf  Acknowledgements.}  In addition to the funding sources mentioned
in the first footnote, we like to thank several institutions
for their hospitality while work on this project was being done: the
Isaac Newton Institute (Weinstein, Xu); the Nankai Institute for
Mathematics (Liu, Weinstein, Xu); IHES, Max-Planck Institut and
Peking University (Xu).  Thanks go
also to Yvette Kosmann-Schwarzbach, Jiang-hua Lu, and Kirill Mackenzie
for their helpful comments.  

\section{Dirac structures on a  Lie bialgebroid}
\label{sec:dirac}

The notion of Lie bialgebroids is a  natural generalization of
that of Lie bialgebras. Roughly speaking, a Lie bialgebroid
is a pair of Lie algebroids ($A$, $A^*$) satisfying
a certain compatibility condition.
  Such a condition, providing
a definition of {\em Lie bialgebroid}, was given in
\cite{MackenzieX:1994}.  We quote here an equivalent formulation from
\cite{K-S:1994}.
\begin{defi}
A {\em Lie bialgebroid} is a dual pair $(A,A^*)$ of vector bundles equipped with
Lie algebroid structures such that the differential $d_*$ on $\Gamma(\wedge^*A)$
coming from the structure on $A^*$ is a derivation of the Schouten-type bracket
on $\Gamma(\wedge^*A)$ obtained by extension of the structure on $A$. 
Equivalently, $d_{*}$ is a derivation for sections of $A$, i.e.,
\begin{equation}
d_{*}[X, Y]=[d_{*}X, Y]+[X, d_{*}Y], \ \ \forall X ,Y \in \gm (A).
\end{equation}
\end{defi}
For a Lie bialgebroid $(A, A^{*})$, the base $P$ inherits a natural
Poisson structure:
\begin{equation}
\label{eq:base}
\{f, \ g\}_{P}=<df , \ d_{*}g>, \ \ \forall f, g\in C^{\infty}(P),
\end{equation}
where $d_{*}: C^{\infty}(P)\lon \gm (A)$ and
$d: C^{\infty}(P)\lon \gm ( A^{*})$ are the
usual differential  operators associated to Lie 
algebroids \cite{MackenzieX:1994}. It is easy to check the identity
\begin{equation}
\label{eq:dfg}
[df, \ dg]=d\{f, \ g\}_{P}.
\end{equation}

Given a Lie bialgebroid ($A$,  $A^{*}$)  over the base
$P$,  with anchors $a$ and $a_{*}$ respectively,
let $E$ denote their   vector bundle direct sum:
$E=A\oplus A^{*}$.
On   $E$, there exist  two natural nondegenerate
bilinear forms, one symmetric and another antisymmetric, which are
defined as follows:

\begin{equation}
\label{eq:pairing}
(X_{1}+\xi_{1} , X_{2}+\xi_{2})_{\pm}=\half (\langle \xi_{1},  X_{2}  \rangle
  \pm \langle \xi_{2} ,  X_{1}\rangle ).
\end{equation}

On $\gm (E)$, we  introduce a bracket  by

\begin{equation}
\label{eq:double}
[e_{1}, e_{2}]=([X_{1}, X_{2}]+L_{\xi_{1}}X_{2}-L_{\xi_{2}}X_{1}-d_{*}(e_{1},
 e_{2})_{-})
+ ([\xi_{1} , \xi_{2}]+L_{X_{1}}\xi_{2}-L_{X_{2}}\xi_{1} +d(e_{1}, e_{ 2})_{-}),
\end{equation}
where $e_{1}=X_{1}+\xi_{1}$ and $e_{2}=X_{2}+\xi_{2}$.\\\\

Finally, we let $\rho : E\lon TP$ be the bundle map defined by
$\rho =a +a_{*}$. That is,
\begin{equation}
\rho (X+\xi )=a(X)+a_{*} (\xi  ) , \ \ \forall X\in \gm (A) \mbox{ and }
\xi \in \gm (A^{*}).
\end{equation}

 When $(A, A^{*})$ is a Lie bialgebra $(\frakg , \frakg^{*})$,
the  bracket above reduces to the  famous Lie
bracket  of Manin on the double $\frakg \oplus \frakg^{*}$.
On the other hand, if $A$ is the tangent bundle Lie algebroid
$TM$ and $A^{*}=T^{*}M $  with zero bracket,
then Equation (\ref{eq:double}) takes  the form:
\begin{equation}
\label{eq:courant}
[X_{1}+\xi_{1} , X_{2}+\xi_{2}]=[X_{1} ,X_{2}]
+\{ L_{X_{1}}\xi_{2} -L_{X_{2}}\xi_{1} +d(e_{1}, e_{ 2})_{-} \}.
\end{equation}
This is the bracket first introduced by Courant \cite{Courant:1990}, then
generalized to  the  context of the formal variational calculus by Dorfman
\cite{Dorfman}. In general, $E$ together with this bracket
and the  bundle map $\rho$  satisfies certain properties as outlined in
the following:

\begin{pro} \cite{LWX}
\label{pro:LWX}
Given a Lie bialgebroid $(A, A^{*})$, let $E=A\oplus A^*$.
Then $E$, with
the  nondegenerate  symmetric  bilinear form
 $( \cdot , \cdot )_{+}$,  the  skew-symmetric
bracket $[\cdot , \cdot ]$ on $\gm (E)$
and the bundle map $\rho :E\lon TP$ as introduced above,
satisfies the following    properties:
 \begin{enumerate}
\item For any $e_{1}, e_{2}, e_{3}\in \gm (E)$,
$[[e_{1}, e_{2}], e_{3}]+c.p.=\D  T(e_{1}, e_{2}, e_{3});$
\item  for any $e_{1}, e_{2} \in \gm (E)$,
$\rho [e_{1}, e_{2}]=[\rho e_{1}, \rho  e_{2}];$
\item  for any $e_{1}, e_{2} \in \gm (E)$ and $f\in C^{\infty} (P)$,
$[e_{1}, fe_{2}]=f[e_{1}, e_{2}]+(\rho (e_{1})f)e_{2}-
(e_{1}, e_{2})\D f ;$
\item $\rho \smalcirc \D =0$, i.e.,  for any $f, g\in C^{\infty}(P)$,
$(\D f,  \D  g)=0$;
\item for any $e, h_{1}, h_{2} \in \gm (E)$,
  $\rho (e) (h_{1}, h_{2})=([e , h_{1}]+\D (e ,h_{1}) ,
h_{2})+(h_{1}, [e , h_{2}]+\D  (e ,h_{2}) )$,
\end{enumerate}
where
\begin{equation}
\label{eq:T0}
 T(e_{1}, e_{2}, e_{3})=\third ([e_{1}, e_{2} ], e_{3})+c.p.,
\end{equation}
and
$\D :  C^{\infty}(P)\lon \gm (E)$
is  the map: $\D=d_{*}+d $.
\end{pro}

Objects satisfying the above properties are called {\em Courant algebroids}
in \cite{LWX}. In other words, we have:

\begin{thm}
If  $(A, A^{*})$ is a Lie bialgebroid, then
$E=A \oplus A^{*}$ together with $([\cdot , \cdot ], \rho , (\cdot , \cdot)_{+})$
is a Courant algebroid.
\end{thm}

In this case, $E$ is called the {\em double of the Lie bialgebroid}.
\begin{defi}
Let $E =A\oplus A^*$ be the double of a Lie bialgebroid
$(A, A^{*})$. A subbundle $L$ of $E$ is
called {\em isotropic} if it  is isotropic under the
 symmetric bilinear form $( \cdot , \cdot )_{+}$. It is called  {\em integrable}
if $\gm (L)$ is closed under the bracket $[\cdot , \cdot ]$.
A {\em Dirac structure}, or {\em Dirac subbundle}, of
the Lie bialgebroid $(A, A^{*})$ is  a subbundle $L\subset E$
which is maximally isotropic and integrable.
\end{defi}

The following proposition follows immediately from the definition of
Dirac structures, and
Properties (i)-(v) of Proposition \ref{pro:LWX}.

\begin{pro}
\label{prop}
Suppose that $L$ is an integrable isotropic subbundle of
a Courant algebroid \\ $(E, \rho , [\cdot , \cdot ], ( \cdot , \cdot
)) $.
Then $(L, \rho |_{L} , [\cdot , \cdot ])$
is a Lie algebroid.
In particular,  any Dirac subbundle itself is a Lie algebroid.
\end{pro}

\section{Dirac structures and Poisson reduction}
Suppose  that  $(A, A^{*})$ is a  Lie bialgebroid over the base
manifold $P$,  with anchors $a$ and $a_{*}$ respectively.
Let  $E=A\oplus A^{*}$ denote its double, and
let  $L\subset E$ be a
Dirac subbundle.    Clearly, $L\cap  A$ is  a (singular) 
subalgebroid of $A$, and therefore $\cald=a(L \cap A)$ is an integrable
(singular) distribution on $P$. We call $\cald$ the characteristic
distribution of $L$.  Let $\calf$ denote the corresponding (generally
singular) foliation of $P$.
\begin{defi}
A Dirac subbundle $L \subset E$ is called \regular if  its characteristic
distribution $\cald$  induces a simple 
foliation. Here by a simple foliation, we mean a 
 regular foliation $\calf$ such that
 $P/\calf $ is a nice manifold such that
the projection is a submersion.

A function $f \in C^{\infty}(P)$ is called $L$-admissible if there is
a  section in $\gm (A)$, denoted by $Y_f$ ($Y_{f}$  may not be unique),
 such that $ Y_f + df \in \gm(L)$.
We write $C^{\infty}_{L}(P)$ for the set of all $L$-admissible functions.
\end{defi}

Let $L\subset E$ be a \regular Dirac structure.
Then $f$ is $L$-admissible iff $ f$ is
constant along $\calf$, i.e.,
$$ C^{\infty}_{L}(P) \cong C^{\infty}(P/ \calf ). $$

For any $f, \ g \in  C^{\infty}_{L}(P)$ admissible, define a
bracket by
\begin{equation}
\label{eq:reduced}
      \{f,g\}=\rho (e_f)g,
\end{equation}
where $e_f=Y_f+df$, which is unique up to
a section of $L \cap A$.

\begin{thm}
\label{thm:main01}
Suppose that $L$ is a \regular Dirac structure. The bracket (\ref{eq:reduced})
defines a Poisson structure on $C^{\infty}(P/ \calf )$.
\end{thm}
\pf  From the definition, we have
$$ \{f,g\}=<Y_{f},dg> +<df,d_{*}g>,  \\\ \forall f,g \in C^{\infty}_{L}(P).$$

On the right  hand  side, the first term is skew-symmetric since
$L$ is isotropic.  The second term
is just the Poisson bracket on  $C^{\infty} (P)$  as defined by
Equation (\ref{eq:base}).  Hence, $\{\cdot,\cdot\}$ is
skew-symmtric.

Next, we prove  that $C^{\infty}_{L}(P)$ is
closed under this bracket and   the Jacobi identity holds.

Let $[e_f,e_g]^*$ denote
the component of $[e_f,e_g]$ on $ \gm (A^*)$.
According to Equation (\ref{eq:double}),
\be
     [e_{f}, e_{g}]^{*}&=&
         [df , dg]+L_{Y_{f}}dg-L_{Y_{g}}df +d(e_{f}, e_{ g})_{-}\\
&=&
         [df , dg]+L_{Y_{f}}dg-L_{Y_{g}}df +d(e_{f}, e_{ g})_{-}
+d(e_{f}, e_{ g})_{+}\\
    & =&   d\{f,g\}_{P}+ d<Y_f, dg>\\
       & =&d\{f,g\}.
\ee
This means that $\{f,g\}$ is also $L$-admissible and one can take
$[e_f,e_g]$ as $e_{\{f,g\}}$. It follows that
$$
\{\{f,g\},h\}= \rho(e_{\{f,g\}})h\\
         =\rho( [e_f,e_g])h\\
         =[\rho(e_f),\rho(e_g)]h\\
         =\{f,\{g,h\}\}-\{g,\{f,h\}\}.
$$
That is, $\{\cdot,\cdot\}$ defines a  Poisson structure on $P/\calf $. \qed

In what follows, we apply the result above to a special class
of Lie bialgebroids: Lie bialgebroids of Poisson manifolds.
Moreover, we will prove that in this case Dirac structures, roughly speaking,
are in one-one correspondence with Poisson structures
on quotient spaces of the Poisson manifold.

Given a Poisson manifold $(P, \pi )$, its cotangent bundle $T^* P$ 
inherits  a natural Lie algebroid structure, called
the cotangent Lie algebroid of the Poisson manifold $P$ \cite{CDW}.
On the other hand, the tangent bundle $TP$  is a Lie algebroid
in an evident sense. It is known 
that they constitute a Lie bialgebroid \cite{MackenzieX:1994}.
 For simplicity,
we will  use 
 $(TP, T^{*}P; \pi )$ to denote this Lie bialgebroid.

 As a special case of Theorem \ref{thm:main01}, therefore,
any \regular Dirac structure  $L$ on
its double $TP\oplus T^{*}P$ will induce a  Poisson
structure on the quotient space $P/\calf$.

Conversely, suppose that
$\cald$ is an integrable distribution  on $P$
with foliation $\calf$, which is simple.
Assume that 
$M=P/\calf$ has  a Poisson  structure.
 Let $    J: P \lon M $ denote  the  natural projection.

To keep under control the fact that $J$ is not a Poisson map, we define a ``difference'' bracket $C^{\infty}(M)\times C^{\infty}(M)\lon C^{\infty}(P)$ by:
\begin{equation}
\label{eq:b1}
\{f,g\}_{1}=J^{*}\{f,g\}-\{J^{*}f,J^{*}g\}_{P}, \ \ \
\forall  f, \ g \in C^{\infty}(M).
\end{equation}

It is easy to see, 
 by using the Leibniz identity of Poisson brackets,
 that this  bracket defines a   skew-symmetric bilinear
form on the conormal bundle $\cald^{\perp }$, which in turn
induces a   bundle map
 \begin{equation}
\label{eq:lambda}
 \Lambda :\cald^{\perp}\lon TP/\cald.
\end{equation}
Let $pr: TP \lon TP/\cald$ be  the natural projection.

Define a  subbundle $L\subset TP\oplus T^{*}P$ by
\begin{equation}
\label{eq:L}
 L=\{(\nu , \ \xi )|pr (\nu )= \Lambda \xi ,  \ \forall
 \nu \in TP , \xi \in \cald^{\perp} \}.
\end{equation}

It is clear that $L$ is a maximal isotropic subbundle of $TP\oplus T^{*}P$,
 and $ C^{\infty}_{L}(P) \cong C^{\infty}(M)$.
For any $f\in C^{\infty}_{L}(P)$, it is easy to see, from
definition, that there exists a vector field $Y_{f}\in \calx (P)$
such that $Y_{f}+df\in \gm (L)$. And in fact, $\gm (L) $ is
spanned by all those sections  of the form  $g (Y_{f} +df)$, for $f\in
C^{\infty}_{L}(P)$ and $g\in C^{\infty}(P)$. To prove that $L$ is
integrable, it suffices to show that the bracket is closed
for those sections having  the form
 $ Y_{f} +df$  according to  Property (iii) of Proposition \ref{pro:LWX},
since $L$ is isotropic.

Given any $f $ and
$g \in C^{\infty}_{L}(P)$. Let   $e_{f} =Y_{f}+df $ be
a section  in $ \gm (L)$. Similarly, let $e_{g}=Y_{g}+dg$ and
$e_{\{f , g\}}=Y_{\{f, g\}}+d\{f, g\}\in \gm (L)$.
It is easy to check that
$$\rho(e_f)g=\{f,g\}.$$
By virtue of the Jacobi identity, we have

\begin{equation}
\label{eq:fgh}
\rho([e_{f},e_{g}]-e_{\{f,g\}})h=0,
\\\        \forall f,g,h \in  C^{\infty}_{L}(P).
\end{equation}
Since the component of $[e_{f},e_{g}]$ on $\gm (T^{*}P)$ is $d\{f, g\}$
according to the proof of Theorem \ref{thm:main01},
$[e_{f},e_{g}]-e_{\{f,g\}} $ is a section in $\gm (TP)$.
Thus Equation (\ref{eq:fgh})  implies that
$$[e_{f},e_{g}]-e_{\{f,g\}} \in \gm (L \cap TP) \subset \gm (L). $$
Hence, $ [e_{f},e_{g}] \in  \gm (L)$,
 and therefore $L$ is integrable.

This proves the following
\begin{thm}
\label{thm:main02}
Suppose that $P$ is a Poisson manifold. There is a one-one
correspondence between \regular Dirac structures in
the double  $E=TP \oplus T^{*}P$ and
 Poisson structures on a  quotient  space  $P/\calf$.
\end{thm}
{\bf Remark} (1). This is a generalization of  a result of Courant
\cite{Courant:1990}. As a main motivation for the  introduction 
of Dirac structures, Courant proved that
a Dirac structure on $TP\oplus T^{*}P$, when $P$ is equipped with
the zero Poisson structure, induces a Poisson bracket on a quotient
space.
In fact, the first part of 
Theorem \ref{thm:main02} can be obtained by  reducing 
to the zero Poisson case. This can be seen as follows.
If  $P$ has a non-trivial Poisson structure, the double
$TP\oplus T^{*}P$, as a Courant algebroid,
is still  isomorphic to the double studied by Courant.
As a consequence, any Dirac structure will thus  induce
a Poisson structure on the quotient. However,
the converse seems to be   new.

(2). If the Poisson structure on the
quotient $P/\calf$ is induced from that on $P$, i.e.,
the projection $J$ is a Poisson map, then it is simple
to see that   $L=\cald \oplus \cald^{\perp}$. This is called a
{\em null Dirac structure} (see \cite{LWX}).
Thus, we have proved: a foliation $\calf$ on a Poisson
manifold $P$ is compatible with the Poisson structure
iff $L=\cald \oplus \cald^{\perp}$ is a Dirac subbundle
of $E=TP\oplus T^*P$.

(3). The result is even interesting when $\cald= 0$. In this case, a Poisson
structure on the quotient   is simply
 another Poisson structure $\pi_{1}$ on $P$.
  Then, $L$ is   simply the graph of the bundle map
$T^{*}P \lon TP$ induced
from the  bivector field $\pi_{1}-\pi $.


\section{Invariant Dirac structures}

This section is devoted to the study of invariant
Dirac structures of a Poisson group. As an application, we will  give
a new proof for the Drinfel'd theorem on homogeneous spaces.

\begin{lem}
\label{lem:invariant}
Let $P$ be a Poisson manifold with a Lie group $G$-action:
$\{\phi_{k}\}_{ k \in G}$. Write $\phi_{k}(x)=kx$ for $x \in P$.
Suppose that $L \subset TP \oplus T^*P$ is  a \regular Dirac
subbundle.
 Then, $L$ is $G$-invariant  iff both  the characteristic distribution
$\cald$ and the difference bracket $\{\cdot,\cdot\}_1$  defined by Equation
 (\ref{eq:b1}) are
$G$-invariant.
\end{lem}
\pf
Firstly, Suppose that $L$ is $G$-invariant. Then,  $\cald =L\cap TP$
is clearly $G$-invariant.
For any  $k \in G$ and $Y_{f} +df \in \gm (L)$,

$$({\kkk}^{-1}+\kk^{*})(Y_{f}+df)=\kkk^{-1}(Y_{f})+d(\kk^{*}f) \in \gm (L).$$
This means that $\kk^{*}f$ is also $L$-admissible  and
we may take  $ Y_{\kk^{*}f}= \kkk^{-1}(Y_{f})$.
Hence,
\be
   \{\kk^{*}f,\kk^{*}g\}_1(x)&=&<Y_{\kk^{*}f}(x),d(\kk^{*}g)(x)>\\
              &=&<Y_{\kk^{*}f}(x),\kk^{*}(dg(kx))>\\
           &=&<Y_{f}(kx),dg(kx)>\\
           &=&\{f, g\}_{1}(kx)\\
           &=&\kk^{*}\{f,g\}_1(x).
\ee

Conversely, from the assumption, it is easy to see that the bundle map
$\Lambda$ as defined by Equation (\ref{eq:lambda}) is $G$-invariant.
Thus $L$ is $G$-invariant according to Equation (\ref{eq:L}).  \qed
{\bf Remark.}\ We note that in general the group action
does not preserve the bracket on the double
$TP\oplus T^*P$ unless it preserves the Poisson
structure on $P$. However,  as we shall see below, in most interesting cases,
we need to study  a Poisson group action, which does
not preserve the Poisson structure.

\begin{thm}
\label{thm:action}
With the notation above, suppose  that $G$ is a Poisson group and $P$ is a
Poisson $G$-space. Then the following statements are equivalent:
\begin{enumerate}
\item L is $G$-invariant.
\item The $G$-action can be reduced to the quotient space $P/{\cal F}$ such
that the reduced action is also a  Poisson action.
\end{enumerate}
\end{thm}
\pf
By  definition, $\{\cdot,\cdot\}_1$ is $G$-invariant iff
\begin{equation}
\label{eq:p2}
\kk^*\{f,g\}_{P}- \{\kk^*f,\kk^*g\}_{P}= \kk^*\{f,g\} - \{\kk^*f,\kk^*g\}, \ \ \
 \forall f , g\in C^{\infty}(P/\calf ).
\end{equation}

Recall that, for a Poisson Lie group $G$,  a Poisson manifold $P$
with a $G$-action
is a  Poisson $G$-space iff the following equality:
\begin{equation}
\label{eq:p3}
    \kk^*\{f,g\}_{P}(x) - \{\kk^*f,\kk^*g\}_{P}(x)=\{f_x,g_x\}_G(k)
\end{equation}
holds for all $k \in G, x \in P$, and $f,g \in C^{\infty}(P),$
where  the function $f_x \in C^{\infty}(G)$ is  defined by $f_x(k)=f(kx)$.

When $L$ is $G$-invariant, its  characteristic foliation
${\cal F}$  is also $G$-invariant. Hence,  the action can be reduced
to  $P/{\cal F}$. Moreover,
combining  Equations (\ref{eq:p2})
and (\ref{eq:p3}), we get
\begin{equation}
\label{eq:com}
    \kk^*\{f,g\}(x) - \{\kk^*f,\kk^*g\}(x)=\{f_x,g_x\}_G(k),
\forall k \in G, x \in P, \mbox{  and } f, g \in C^{\infty}_L(P).
\end{equation}
This means that $P/{\cal F}$ is a  Poisson $G$-space.

Conversely, if both $P$ and $P/\calf$ are Poisson $G$-spaces,
then  Equations (\ref{eq:p3}) and (\ref{eq:com})
imply  Equation (\ref{eq:p2}), which is equivalent
to $L$ being $G$-invariant.

\qed

Now, by means of the theory of Dirac structures, we are in a position
to explain the Drinfel'd's more or less mysterious theorem regarding
Poisson homogeneous spaces, outlined in his short paper
\cite{dr:poisson} (also see \cite{Lu} for the intepreation of the
associated Dirac structures of Poisson homogeneous spaces in terms of
Lie algebroids).

To  a Poisson Lie group $(G, \pi_{G})$,
there  are associated two Lie bialgebroids.
One is  $(TG,T^*G; \pi_G)$  
with the canonical Lie bialgebroid structure induced
from the Poisson structure on $G$,
 and  the other  is its tangent  Lie bialgebra $(\frakg,\ \frakg^*)$.
Identifying $TG\oplus T^*G$ with the trivial vector bundle
$G\times (\frakg \oplus  \frakg^* )$ by left translations,
it is clear that  there is a  1-1 correspondence between
  maximal isotropic subspaces $L$ of $\frakg\oplus \frakg^*$ and  left
invariant  maximal isotropic subbundles $\bar L$ of $TG\oplus T^*G$.
Moreover, we have

\begin{lem}
\label{lem:bar}
  Given   any $e_{1}, \ e_{2}\in \frakg \oplus \frakg^{*}$,
  let $\bar{e_{1}} $ and $\bar{e_{2}}$ 
 be their corresponding
left invariant sections  in $TG\oplus T^{*}G$. 
Then,
$$[\bar e_1, \ \bar e_2] =[e_1,e_2]^- . $$
\end{lem}
\pf Assume that   $e_{1}=X+\xi$,   $e_{2}=Y+\eta \in \frakg \oplus
\frakg^{*}$, and   $\bar{e_1}=\bar X+\bar{\xi}$,
 $\bar{e_{2}} =\bar Y +\bar{\eta }$.
Here $\bar X, \  \bar Y$ and $\bar \xi, \ \bar \eta$ denote
the corresponding left invariant vector fields and
1-forms on $G$ for  $X,Y \in \frakg$ and $\xi,\eta \in \frakg^*$ respectively.

 From the fact that $[\bar X,\bar Y]=[X,Y]^-$ and
$[\bar \xi, \bar \eta]=[\xi,\eta]^-$
(see \cite{we:dressing}),
it follows that
$$
L_{\bar X}\bar \xi =(ad^*_X \xi)^-, L_{\bar \xi}\bar X=(ad^*_{\xi}X)^- .
$$

This means that
\be
[\bar X, \ \bar \xi]&=&L_{\bar X}\bar \xi -L_{\bar \xi}\bar X
                       + \frac{1}{2}(d^*-d)<\bar X,\bar \xi>\\
           &=&(ad^*_X \xi)^- -(ad^*_{\xi}X)^-\\
           &=&[X,\xi]^-.
\ee

Therefore,
\be
[\bar e_1, \ \bar e_2]&=&[\bar X, \bar Y]+[\bar X,\bar \eta]+[\bar \xi,\bar Y]+
[\bar \xi, \bar \eta]\\
&=&[X,Y]^- +[X,\eta]^- +[\xi,Y]^- +[\xi,\eta]^-\\
&=&[e_1,e_2]^- .
\ee
\qed

An immediate consequence is the following:

\begin{cor}
As above, $L$ is a subalgebra iff $\bar L$ is integrable. That is,
Dirac structures of $\frakg\oplus \frakg^*$ are in
 1-1 correspondence  with  left invariant Dirac
structures of $TG\oplus T^*G$ .
\end{cor}
\pf    Sections  of $\bar L$ are  spanned by those  of the
form $f\bar e$, where $f \in C^{\infty}(G)$.
The conclusion thus follows  immediately from
Lemma \ref{lem:bar} by using
Property (iii) of Proposition \ref{pro:LWX}. \qed 

 Given a left invariant 
 $TG\oplus T^*G$ Dirac stucture $\bar L$ (for the bialgebroid
 associated with the Poisson structure on $G$), its
 characteristic distribution $\bar L \cap TG$ is just the left
 translation of the subalgebra $\frakh=L \cap \frakg$ of $\frakg$. So
 its quotient space is $G/H$, where $H$ is the connected subgroup of
 $G$ with Lie algebra $\frakh$. On the other hand, it is well known
 that $G/H$ is a nice manifold such that the projection is a
 submersion iff $H$ is closed.  In this situation, we call $L$ a {\em
 \closed} $(\frakg,\frakg^*)$ Dirac structure. 
In other words, $L$ is regular if
the left translation of $\frakh=L \cap \frakg$
defines a simple foliation on $G$.
  It is simple to see that $L$ is \closed
 iff $\bar L$ is \regular.  Thus, by Theorem \ref{thm:action}, we
 obtain the following:

\begin{thm}
 \Closed Dirac structures of $\frakg\oplus \frakg^*$ are in
 1-1 correspondence with Poisson homogeneous spaces $G/H$, where $H$
is a  connected closed subgroup of $G$.
\end{thm}

Every homogeneous $G$ space $X$ is of  the form $X=G/\tilde H$,
where $\tilde H $ is a closed subgroup of $G$,
 with $H$ being its connected component at the unit.
 Then $D=\tilde{H}/H$ is a discrete group, and  the
 projection $p: \ G/H\lon X$ is a covering map with structure
 group $D$. 
Any Poisson structure on $X$ can be
pulled back to  $G/H$ such that  $p$   is a   Poisson map.
Moreover, if one is a Poisson  homogeneous space,
so is the other.
 Thus, any Poisson homogeneous $G$-space
$X=G/\tilde H$
induces a \closed  Dirac structure $L$ in $\frakg\oplus \frakg^*$.
It is easy to see that  $L$ is
$Ad_{\tilde H}$-invariant.

Conversely, given a  \closed Dirac structure $L$ of $\frakg\oplus \frakg^*$,
let $G/H$ be its  corresponding Poisson homogeneous $G$-space. Then,
the Poisson structure on $G/H$ can be reduced to a
 homogeneous $G$-space $X=G/\tilde H =(G/H)/D$ iff $L$ is $Ad_D$-invariant,
or equivalently,  is $Ad_{\tilde H}$-invariant.

Thus, we obtain the following:
\begin{thm} \cite{dr:poisson}
 Poisson homogeneous $G$-spaces bijectively
correspond to  pairs $(L,K)$, where  $L$  is a \closed Dirac structure of
$\frakg\oplus \frakg^*$ and $K$ is a closed subgroup of $G$ with
Lie algebra $L \cap \frakg$ such that $L$ is invariant under the
(adjoint,coadjoint) action of $K$.  
\end{thm}

\section{Pullback of Dirac structures}

 This   section is devoted to a discussion of
 pullbacks of Dirac structures.
It will be  used in Section 6 to
extend $(A,A^* )$ Dirac structures on $P$ to ``left-invariant''
$(TG,T^*G)$ 
Dirac structures on the Poisson groupoid $G$.  These
invariant structures are then related to Poisson structures on
quotients of $G$.

Given  two  vector spaces $U$ , $V$,  and  a surjective
linear map $\Phi: U \lon V$,  its dual
  $ \Phi^* : V^* \lon U^*$ is injective  and
 $ \Phi^* (V^*) = (ker\Phi )^{\perp}$.
Write
 $$ \bar \Phi  = \Phi  \oplus (\Phi^*)^{-1} : U\oplus  (ker\Phi )^{\perp}  \lon
 V\oplus  V^*.$$
Clearly, $\bar \Phi$ is a surjective linear map.
Given any  maximal isotropic subspace
$L \subset  V\oplus V^*$, we denote by $ \bar L$
the inverse image $ \bar \Phi^{-1}(L)$.
 Then,
 $\bar L$ is a maximal isotropic subspace of $ U \oplus  U^*$, which
is called the pullback of $L$.

 Similarly  suppose that $A\lon P$ and $B\lon Q$ are vector bundles, and
 $\Phi : A\lon B$ is a surjective  bundle map covering a map $P\to Q$.
Then,  given  any maximal isotropic subbundle $L\subset B\oplus B^*$,
we may define its pullback $\bar L \subset A\oplus A^*$ as
the fiberwise pullback. Then, $\bar L$ is a
 maximal isotropic subbundle, and the restriction
of $\bar \Phi$ on $\bar L$, denoted by the same symbol
$\bar \Phi$, is a vector bundle morphism $\bar L\lon L$.

Given a vector bundle morphsim $\rho :E_{1}\lon E_{2}$,
a section $\bar X$ of $E_{1}$ will be called {\em admissible} if $\rho \bar{X}$
corresponds to a section $X$ in $E_{2}$. In this case,
$\bar X$ and $X$ are said to be  $\rho$-related. If $\rho$ is
surjective, then $\gm (E_{1})$ is spanned  
over $\reals$
(possibly infinite sum but  locally finite)
 by all sections of the
form $f\bar X$ for $f\in C^{\infty}(P)$ and $\bar X\in \gm (E_{1}) $
admissible,  by the partition of unity.

Applying the obervation above to the bundle morphim $\bar \Phi :
\bar L\lon L$,  it follows that
\begin{equation}
\label{eq:pull}
 \gm(\bar L) =span\{f\bar e| \forall \mbox{ admissible } \bar{e} \in \gm (\bar L)
\ \mbox{ and }   f \in C^{\infty}(P)\}.
\end{equation}

\begin{thm}
\label{thm:main3}
 Let $(A, A^*),  (B,B^*)$ be two Lie bialgebroids, and $ \Phi  : A \lon B$
     a  surjective bundle map. Thus the bundle map
 $ \bar \Phi  = \Phi  \oplus (\Phi^*)^{-1} : A\oplus  (ker\Phi )^{\perp}  \lon
 B\oplus  B^*$ is     surjective. 

Given any maximal isotropic subbundle $L\subset  B\oplus B^*$, its pull back
$\bar L={\bar \Phi }^{-1}(L)$ is a Dirac structure of $A\oplus A^*$
iff  $L$ is Dirac structure.
     Moreover, in this case,
 $\bar \Phi: \bar L \lon L$ is a  Lie algebroid morphism.
\end{thm}
\pf Since $\gm (\bar L)$ is spanned by sections of the form $f\bar e$,
it suffices to prove the  following identity:
\begin{equation}
\label{eq:e1}
\bar \Phi [\bar e_1,\bar e_2]=[e_1,e_2] , \ \ \ \ \forall
 \ \mbox{ admissible }\bar  e_1, \ \bar e_2 \in \gm(\bar L),
\end{equation}
according to  Property (iii) of  Proposition \ref{pro:LWX}.

Write $\bar{e_{1}}=\bar{X}+\bar{\xi}$ and $\bar{e_{2}}=
\bar{Y}+ \bar{\eta} $, where $\bar X$ and $\bar Y$ are admissible
sections of $A$ under the map $\Phi$, and $\bar \xi$ and $\bar \eta $
are admissible sections of $(ker\Phi )^{\perp}$ under the map $(\Phi^* )^{-1}$.
Denote by  $X$,  $Y $, and $\xi, \ \eta$ their  corresponding sections
in $B$ and $B^{*}$, respectively.

Since  $\Phi$ is a Lie algebroid morphism, by definition (see \cite{HM:1990}),
\begin{equation}
\label{eq:e2}
 \Phi [\bar X, \bar Y]=[X,Y].
\end{equation}

Moreover, $\Phi$ is a Poisson map, where $A$ and $B$ are equipped with the
Lie-Poisson structures corresponding to the  Lie algebroids $A^*$ and $B^*$
respectively, since $\Phi$ is a Lie bialgebroid morphism.
Thus,
$$\Phi^* \{l_{\xi}, l_{\eta}\}=\{\Phi^* l_{\xi},\ \Phi^*  l_{\eta}\},
$$
where
$l_{\xi}$ and $l_{\eta }$ are
the  linear functions on $B$ corresponding to $\xi, \ \eta \in \gm (B^{*})$.
Therefore,
$$\Phi^{*}l_{[\xi , \eta ]}=\{l_{\bar \xi }, l_{\bar \eta }\}=l_{[\bar{\xi} ,
\bar{\eta }]}. $$
Thus it follows that
\begin{equation}
\label{eq:e3}
(\Phi^*)^{-1}[\bar \xi,\bar \eta]=[ \xi,\eta].
\end{equation}
 That is,
$(\Phi^*)^{-1}: (ker \Phi)^{\perp} \lon B^*$ is also a Lie algebroid
morphism, where $(ker \Phi)^{\perp}$ is considered as a subalgebroid of $A^*$.
Now
$$
[\bar e_1,\bar e_2]=[\bar X, \bar Y]+[\bar X,\bar \eta]+[\bar \xi,\bar Y]+
[\bar \xi, \bar \eta].
$$
According to Equation \ref{eq:double},
$$
[\bar X,\bar \eta]=L_{\bar X}\bar \eta -L_{\bar \eta}\bar X
+ \frac{1}{2}(d^*-d)<\bar X,\bar \eta>
=i_{\bar X}d\bar \eta -i_{\bar \eta}d_*\bar X
- \frac{1}{2}(d_*-d)<\bar X,\bar \eta>.
$$
Since $\Phi$ is a Lie algebroid morphism, $ d\bar \eta$ is admissible
and is  $(\Phi^*)^{-1}$-related to
$d\eta$. Similarly, $d_*\bar X$ is $\Phi$-related to $d_*X$.
Finally, note that $
<\bar{X} ,\bar{\eta}>
=\phi^{*}<X, \eta>$, so
$d_{*}
<\bar{X} ,\bar{\eta}>$  and $d_{*}<X, \eta >$ are $\Phi$-related
while $d
<\bar{X} ,\bar{\eta}>$ and $d<X, \eta >$ are
 $(\Phi^*)^{-1}$-related.  Hence,

$$\bar \Phi [\bar X,\bar \eta]=[X,\eta]. $$
Similarly, $ \bar{\Phi } [\bar \xi , \bar Y]=[\xi , Y]$.
Hence, Equality (\ref{eq:e1})  follows. This concludes the proof of the theorem.
\qed

\begin{ex}
Recall that  a hamiltonian operator on  a Lie bialgebroid $(B, B^* )$ 
 is a   skew-symmetric
two-form in $\gm (\wedge^{2}B^* )$
 satisfying  the following Maurer-Cartan type equation \cite{LWX}:
$$dI+\half [I, I]=0.$$
In particular, $I$ is called a strong hamiltonian operator  if both $dI$ and $[I
, I]$ vanish.
Given a  two-form $I\in \gm (\wedge^{2}B^* )$,  the
graph  $L_{I}=\{ X+I^{b} X|X\in  B \}$ of  its  induced bundle map    is a Dirac
 structure
iff $I$ is a  hamiltonian operator.

Suppose that $\Phi: A\lon B$ is  a Lie bialgebroid morphism and $I\in 
\gm (\wedge ^{2}B^* )$ a (strong) hamiltonian operator.
 $\Phi^*$ pulls a   two-form
in $\gm (\wedge^{2}B^*)$  back to a   two-form in  $\gm ( \wedge^{2}A^* ) $.
It is easy to see that  $\Phi^*I$ is then  a (strong) hamiltonian operator in $A
$.
Moreover, the pull-back of the corresponding  Dirac structure  $L_{I}$ 
is exactly  the graph $L_{\bar I}$  of the (strong) hamiltonian operator
 $ {\bar I}=\Phi^*I$.
\end{ex}

\begin{ex}
Given a surjective submersion $\phi: M\lon N$ of 
 manifolds $M$, $N$, its  derivative 
defines a Lie  algebroid morphism $\Phi : TM\lon TN$.
This  is also a Lie bialgebroid morphism between  $(TM, T^{*}M)$
and  $(TN, T^{*}N)$, where $T^{*}M$ and $T^{*}N$ 
are considered  as the cotangent algebroids for the zero Poisson
 structure.
According to Courant \cite{Courant:1990},  in this case,
 a Dirac structure simply  corresponds to a foliation on the manifold
together with a family of closed 2-forms on the leaves. The pullback 
of the  Dirac structure
 is just the pullback  of the foliation by $\phi$ together
with the pullback of two-forms.
\end{ex}
{\bf Remark.}  A    Dirac structure on a vector space is 
equivalent to a two-form on a subspace \cite{Courant:1990}. 
Thus, 
we   could  also  pull back an isotropic
 subbundle by any  bundle map, not
just a surjection.
Of course, the pullback might not be continuous if  the map
is not of constant  rank.   Moreover, it is not clear whether the integrability
condition is preserved.

\section{Left invariant Dirac structures on Poisson groupoids}
To generalize Drinfel'd's theorem on homogeneous spaces from Poisson
groups to Poisson groupoids, we will first extend the notion of left
invariant Dirac structure from Poisson groups to groupoids.

Let $(\poidd{G}{P}; \alpha,\beta)$
 be a  Poisson groupoid, with Lie algebroid $A$.
Here  $\gm (A) $ is identified with  
 left-invariant vector fields on the groupoid.
The  dual bundle  $A^*$  can be naturally
identified with the conormal bundle of the
identity space $P$ in the groupoid, and therefore
inherits a Lie algebroid structure according
to Weinstein \cite{we:coisotropic}.
Moreover, it was shown in \cite{MackenzieX:1994}
that   $(A, A^*)$ is a Lie bialgebroid, which
is called the  tangent Lie bialgebroid of $G$.

By $T^{\alpha}G$, we denote the subbundle of $TG$ consisting
of all vectors tangent to $\alpha$-fibers.
The group  $B(G )$ of bisections of  $G$ (submanifolds which
 project diffeomorphically to $P$ by $\alpha$ and $\beta$)
acts naturally  on $G$ by left multiplication:
 $$l_{K}x=K\cdot x, \ \ \ \forall K\in B(G ), \ x\in G. $$
As usual,  the action lifts naturally
 to actions on   $TG$ and $T^{*}G$ which leave
$T^{\alpha}G$ and $(T^{\alpha }G)^{\perp}$ invariant.

Define a   map 
$\Phi : T^{*}G \lon A^{*}$ as follows:  
given any $x\in G$  and  $\xi \in T_{x}^{*}G $,
set $\Phi (\xi )\in A_{p}^{*}$ with $p=\beta (x)$ such that
\begin{equation}
\label{eq:Phi}
<\Phi \xi , v>=<\xi , Tl_{x}v> , \ \ \ \forall v\in A_{p}.
\end{equation}
Then, 
\begin{equation}                         
\matrix{&\Phi&\cr
        T^* G &\vlra&A^*\cr
        &&\cr
        \Bigg\downarrow&&\Bigg\downarrow\cr
        &&\cr
        G&\vlra&P\cr
        &\beta &\cr}
\end{equation}
is a bundle map.

\noindent{\bf Remark.} In terms of symplectic groupoids \cite{CDW},
$\Phi$ is just the $\beta$-map of
the cotangent groupoid $\poidd{T^{*}G }{A^*}$ \cite{CDW}.
However, our  symplectic structure on $T^*G$ differs
a minus sign from the one on the symplectic
groupoid 
 $\poidd{T^{*}G }{A^*}$. Therefore,  $\Phi$ is
a Poisson map.

Another interesting way to think of $\Phi$ is as the
momentum map for the lifted right action of $B(G )$ on
the cotangent bundle $T^{*}G$. Here $A^{*}$ is
considered as a subset of $\gm (A)^{*}$ in the form of delta-distributions.
 It can be checked that  the
 image of the momentum map  is just $A^*$.
This also indicates that $\Phi$ should be  a Poisson map.

The following lemma lists some basic properties of
$\Phi$.

\begin{lem}
\begin{enumerate}
\item $ker \Phi =(T^{\alpha}G)^{\perp}$;
\item for any $X \in \gm (A)$, $\Phi^{*}X$, as a section in
$TG$, is exactly the left invariant tangent vector field
$\bar{X}$ obtained by left translating  $X$ along $\alpha$-fibers;
\item $\Phi$ is $B(G )$-invariant, i.e.,
$$\Phi (l_{K}^{*}\xi )=\Phi (\xi ), \ \ \ \forall K\in B(G ), \mbox{ and }
\xi \in T^{*}G ; $$
\item 
$\Phi$ is a Lie bialgebroid morphism, where $(T^{*}G , TG )$
is equipped with the natural Lie bialgebroid structure
associated with the Poisson structure $G$ while $(A^{*}, A)$
is the flipping of the tangent Lie bialgebroid $(A, A^{*})$
\footnote{In order to be consistent with previous
notation for Lie bialgebroid morphisms, 
we have flipped both Lie bialgebroids here.}.
\end{enumerate}
\end{lem}
\pf The proof for (i)-(iii) is obvious, and is left for the reader.
For (iv),  since $\Phi$ is already known to be a Poisson map,
it suffices to show that it is a Lie algebroid morphism.
This is, however, quite clear since the Lie algebroid
structure on $A^*$ is defined in terms of  the Lie algebroid
structure 
on $T^*G$ by identifying $A^*$ with the
conormal bundle of $P$ in $G$ \cite{we:coisotropic}. In fact,
$T^*G$ is a LA-groupoid in terms of Mackenzie \cite{Mackenzie:private}.  \qed

\begin{pro}
\label{pro:Poi1}
 A Dirac structure $\bar L \subset TG \oplus T^*G$ is the
 pullback of a Dirac structure in  $A\oplus A^*$ iff
\begin{enumerate}
\item  $\bar L$ is $B(G)$-invariant, and
\item  $(T^{\alpha}G)^{\perp} \subset \bar L  $.
\end{enumerate}
\end{pro}
\pf Using $l_{K}$, $\forall K\in B(G )$, we denote the
$B(G )$-action on $TG \oplus T^{*}G $.
 As in the previous  section, 
let $\bar \Phi$ be the map $(\Phi^{*})^{-1}\oplus \Phi :
T^{\alpha}G \oplus T^{*}G \lon A\oplus A^{*}$.
It is obvious that $\bar \Phi$ is also invariant under the $B(G )$-action,
i.e.,
$$\bar{\Phi} \smalcirc l_{K}=\bar{\Phi}, \ \ \ \forall K\in B(G ).$$

Suppose that $\bar L$ is a Dirac structure in $TG \oplus T^* G$.
If $\bar L$  is the pullback of a Dirac structure
$L$ in $A\oplus A^*$. Then, $\bar{L}=\bar{\Phi }^{-1}(L)$.
Clearly, $\bar{L}$ is $B(G )$-invariant since $\bar \Phi$ is invariant.
 Since $(T^{\alpha}G)^{\perp}=ker \Phi $,  it follows that
  $(T^{\alpha}G)^{\perp} \subset \bar L $.

Conversely,  the condition
  $(T^{\alpha}G)^{\perp} \subset \bar L $ implies that
$\bar L\subset
T^{\alpha}G \oplus T^{*}G $,
since $\bar L$ is isotropic.
Since $\bar L$ is $B(G )$-invariant,
it follows that
$\bar{\Phi} (\bar{L}|_{K \cdot x})=
\bar{\Phi} (l_{K} \bar{L}|_{x}) =
\bar{\Phi} (\bar{L}|_{x}) $.  Therefore,
$\bar{\Phi} (\bar{L}|_{x}) $  depends only  on the base point
$p=\beta (x)$, and thus  defines  a subspace $L_{p}$ in
$A_{p}\oplus A^{*}_{p}$, which is easily seen to be maximal isotropic.
Thus we obtain a maximal isotropic subbundle $L$ in $A\oplus A^*$
suct that $\bar{L}={\bar \Phi}^{-1}(L)$. According to theorem \ref{thm:main3},
$L$ must be  a Dirac subbundle. \qed

When $G$ is a Poisson {\em group}, pullback Dirac structures in the
sense above  correspond exactly
 to  left invariant Dirac strucures as discussed  in Section 4. For this
reason, we shall
call any such pullback a {\bf left invariant Dirac structure}.
Suppose that  $\calf $ is  a foliation
on $G$ with distribution $\cald  \subset TG$ such that $G /\calf$ is
a nice manifold.
 According to Theorem \ref{thm:main02}, every Poisson structure on $G/\calf$
corresponds to a Dirac structure $\bar{L}\subset TG \oplus T^* G$.

\begin{thm}
\label{thm:main4}
  $\bar L$ is the pullback of a Dirac structure $L\subset
A\oplus A^{*}$  iff
\begin{enumerate}
\item $\calf$ is $B(G )$-invariant;
\item $\{ \cdot , \cdot \}_{1}$ is $B(G )$-invariant; and
\item $\cald \subset T^{\alpha}G $ and $\{\alpha^{*}f , \ g\}_{1}=0$, $\forall
    f\in C^{\infty}(P)$ and $g\in C^{\infty}
(G/\calf )\cong C^{\infty}_{\bar L}(G )$. 
\end{enumerate}
\end{thm}
\pf Suppose that $ \bar L$ is the pullback of a Dirac structure
$L$ under $\Phi$.
Then $\bar L$ is $B(G )$-invariant. Thus, $\cald=\bar{L}\cap TG$ is
also $B(G )$-invariant. According to Lemma \ref{lem:invariant},
 $\{ \cdot , \cdot \}_{1}$ is $B(G )$-invariant.
For (iii), we note that $\bar L\subset T^{\alpha}G \oplus T^{*}G$
by Proposition \ref{pro:Poi1} since $\bar L$ is isotropic. Thus it follows that
$\cald=\bar{L}\cap TG\subset T^{\alpha}G $.
Also note that $\alpha^{*}f$ is constant along $\alpha$-fibers,
so  $\cald\subset T^{\alpha}G $ implies that $\alpha^* f$ is
admissible. Therefore, $\{\alpha^{*}f, g\}_{1}$ is well-defined.
Now since  $d\alpha^{*}f\in (T^{\alpha}G)^{\perp}=ker\Phi \subset \bar{L}$,
 we may choose  $Y_{\alpha^{*}f}=0$.  Thus,
 $\{\alpha^{*}f , \ g\}_{1}=Y_{\alpha^{*}f}g=0$.

Conversely, assume that (i)-(iii) hold. So, for any $f\in C^{\infty}(P)$,
$\alpha^{*}f$ is $\bar L$-admissible.
As in Section 3, let   $Y_{\alpha^{*}f}\in \calx (G )$ be any vector field
that $Y_{\alpha^{*}f}+d\alpha^{*}f\in \gm (\bar L)$. Then,
$Y_{\alpha^{*}f}g=
 \{\alpha^{*}f , \ g\}_{1}= 0$ for all admissible $g$.
This  implies that
  $Y_{\alpha^{*}f}\in \gm (\cald)\subset \gm  (\bar{L})$. Thus,
$d\alpha^{*}f\in \gm (\bar L)$. Since $\ker\Phi=(T^{\alpha }G)^{\perp}$
is spanned by  all  such vectors, it follows that

$$ker \Phi \subset \bar L. $$
By Lemma \ref{lem:invariant},
(i)-(ii) imply  that $\bar L$ is $B(G )$-invariant. Thus
$\bar L$ is the  pullback of $L$ according to Proposition \ref{pro:Poi1}  . \qed

\section{Poisson actions}
Let $\poidd{G}{P}$ be  a Poisson groupoid with Poisson
tensor $\pi_{G}$. Suppose that $G$ acts on
a Poisson  manifold
   $(X, \pi ) $ equipped with a moment  map $J : X\lon P$.
  Here, the action is a map
$$m:  (g,x)\mapsto g\cdot x$$
 from  $G\times_P X =\{(g,x)\in
G\times X|\beta(g)=J(x)\}$ to $X$, satisfying the usual
condition
$g\cdot(h\cdot x) = (gh)\cdot x$.  The action is a Poisson action if
its graph $\Omega =\{(g, x, g\cdot x)|\beta (g) =J(x) \}$
is a coisotropic submanifold of
$G\times X\times \overline{X}$ \cite{mi-we:moments}.

For example, consider a complete 
 Poisson group  $H$ with dual Poisson group $H^*$.
Then  $G=H H^{*}$ is a 
  symplectic groupoid over   $H^*$. Suppose that $X$ is 
a Poisson $H$-space with an equivariant momentum
map $J: X \lon H^*$. Then,   $(X, J)$ is a Poisson
$G$-space under  the $G$-action \cite{WX}:
$$(g, u) \cdot x = gx,  \ \ \  g \in H, \ \ u   \in H^*
\mbox{  and  } \  x \in X \mbox{ such that  }\ J(x)=  u. $$

Any bisection $K$ of the groupoid  induces   a diffeomorphism of $G$ by
left multiplication by $K$. This is denoted by $l_{K}$.
On the other hand, we also can define a diffeomorphism
of $X$ by $x\lon K\cdot x, \ \forall x\in X$. This is also  denoted
by the same symbol $l_{K}$ when confusion is unlikely.

A section of the moment map is a submanifold of $X$
to which  the restriction of $J$  is a diffeomorphism.
Similarly,  a local section of $J$ is a submanifold
$\caly$ of $X$ to which the restriction of $J$ is
a diffeomorphism between $\caly$ and an open subset $J (\caly )$ of $P$.
If $J$ is a submersion, through any point of $X$ there
always exists a local section of $J$.
A  global  section $\caly$  induces a map, denoted by $r_{\caly}$,
from $G$ to $X$ given by: $g\lon g\caly$.
Similarly, given   any compatible $(g, x)\in G *X $ (i.e. $\beta (g)=J(x)$),
 a local section $\caly$  through the point  $x$
induces   a map $ r_{\caly}$ from
a neighborhood of $g$ to a neighborhood of $x$ in the same way.

The main theorem is the following:

\begin{thm}
\label{thm:main1}
Suppose that $\poidd{G}{P}$ is a Poisson groupoid acting on
a Poisson manifold equipped with a moment map $J: X \lon P$.
This is a Poisson action iff
\begin{enumerate}
\item For any $f\in C^{\infty}(P)$,
\begin{equation}
\label{eq:right}
X_{J^{*}f}(x)=(r_{x})_{*} X_{\alpha^{*}f}(u),
\end{equation}
where $x\in X$ and $u=J(x)$.
\item For any compatible $(g, x)\in G *X$,
\begin{equation}
\label{eq:multi}
 \pi (gx)=(l_{K })_{*} \pi (x) + (r_{\caly })_{*} \pi_{G } (g)
- (r_{\caly })_{*}(l_{K })_{*}\pi_{G}(u),
 \end{equation}
\end{enumerate}
where $u=\beta (g)=J(x)$, $K $ is any local  bisection of
$G$ through $g$, and $\caly$ any local section through
the point $x$.
\end{thm}
\noindent{\bf Remark} (1).\  Here $r_{x}$ denotes the map: $g\lon g\cdot x$
from $\beta^{-1}(u)$ to $X$.
Since $X_{\alpha^{*}f} (u) $ is tangent to $\beta^{-1}(u)$,
the  right hand side of Equation (\ref{eq:right})  makes sense.

(2).\ When $G$ is a Poisson group, i.e. $P$ reduces to a point,
the first condition is satisfied automatically, and
the second one reduces to the usual condition
of Poisson action since $\pi_{G }(u)=0$.

\begin{pro}
\label{pro:1.2}
Under the hypotheses of Theorem  \ref{thm:main1}, if
 the action is a Poisson action,
then $J: X\lon P$ is a Poisson map.
\end{pro}
\pf  Let  $x\in X$, $J(x)=u$ and  $\xi ,\eta \in T_{u}^{*}P $
be any covectors  at $u$.
For any $g\in G$ with $\beta (g)=u$, as covectors at
the point $(g, x, gx)$, both   $(-\beta^{*}\xi, J^{*}\xi, 0)$ and
$(-\beta^{*}\eta , J^{*}\eta , 0)$ are conormal to $\Omega$.
Therefore,
$$\pi_{G }(-\beta^{*}\xi ,-\beta^{*}\eta )+\pi(J^{*}\xi ,J^{*}\eta )=0, $$
which implies that $(J_{*}\pi (x))(\xi,  \eta )=-(\beta_{*}\pi_{G} (g))
(\xi,  \eta )=\pi_{P}(u)(\xi,  \eta )$. That is $J_{*}\pi =\pi_{P}$.
In other words, $J$ is a Poisson map.  \qed

\begin{lem}
\label{lem:1.3}
Let $x\in X $ and  $J (x)=u $.  Suppose that  $\delta_{x} \in T_{x}
X$ and  $\delta_{u}\in T_{u}G $ such   that
$(\delta_{u}, \delta_{x})\in T_{(u,x)} (G *X) $.
Let $\delta_{x}'=m_{*} (\delta_{u}, \delta_{x})\in T_{x}X$.
Then,
\begin{equation}
\label{eq:basic}
\delta_{x}'=\delta_{x}+(r_{x})_{*}(\delta_{u}-\epsilon_{*}J_{*} \delta_{x}),
\end{equation}
where $\epsilon : P\lon G $ is the inclusion  of the unit space.
\end{lem}
\pf As a tangent vector in 
$ T_{(u,x)}  (G *X) $, $(\delta_{u}, \delta_{x})$
can be split into the sum of $(\epsilon_{*}J_{*}\delta_{x}, \delta_{x})$ and
$(\delta_{u}-  \epsilon_{*}J_{*}\delta_{x}, 0)$. It is easy to see that
$m_{*} ( \delta_{u}-  \epsilon_{*}J_{*}\delta_{x}, 0 )=
(r_{x})_{*}(\delta_{u}-\epsilon_{*}J_{*} \delta_{x})$ and
$m_{*}(  \epsilon_{*}J_{*}\delta_{x}, \delta_{x})= \delta_{x}$.
This proves the lemma. \qed

\begin{lem}
Suppose that  $(g, x, z)$ with $z=gx$ is any
point in $\Omega $, and $(\delta_{g}, \delta_{x}, \delta_{z})$  any tangent
vector of  $\Omega$ at this point.
Let $K$ be any  (local) bisection of $G$ through the
point $g$ and $\caly$ any (local) section of $J$
through the point $x$.
Then,
\begin{equation}
\label{eq:diag}
\delta_{z}=r_{\caly_{*}} \delta_{g}+l_{K_{*}}\delta_{x}-
l_{K_{*}} r_{\caly_{*}}\epsilon_{*}J_{*}\delta_{x} .
\end{equation}
\end{lem}
\pf Let $u=\beta (g) =J(x)$, and  $\delta_{u}=l_{K^{-1}_{*}}\delta_{g}$.
Then, $\delta_{u}$ is a  tangent vector of $G$ at $u$.
It is clear that
\be
\delta_{z}&=&m_{*}(\delta_{g}, \delta_{x})\\
&=&l_{K_{*}}m_{*}( l_{K^{-1}_{*}}\delta_{g} , \delta_{x}) \\
&=&l_{K_{*}}m_{*}(\delta_{u}, \delta_{x}) \ \ \mbox{(By Lemma \ref{lem:1.3})} \\
&=&l_{K_{*}}( \delta_{x}+(r_{x})_{*}(\delta_{u}-\epsilon_{*}J_{*} \delta_{x}))\\
&=&l_{K_{*}}( \delta_{x}+r_{\caly_{*} }(\delta_{u}-\epsilon_{*}J_{*} \delta_{x}))\\
&=&l_{K_{*}}( \delta_{x}+r_{\caly_{*} }(l_{K^{-1}_{*}}\delta_{g}
-\epsilon_{*}J_{*} \delta_{x}))\\
&=&r_{\caly_{*}} \delta_{g}+l_{K_{*}}\delta_{x}-
l_{K_{*}} r_{\caly_{*}}\epsilon_{*}J_{*}\delta_{x}.
\ee
\qed

\begin{cor}
Suppose that $(g, x, z)\in \Omega $,
$K$  any  (local) bisection of $G$ through the
point $g$, and $\caly$ any (local) section of $J$
through the point $x$.
For any $\zeta\in T_{z}^{*}X$,
the covector $(-r_{\caly}^{*}\zeta ,\  J^{*}\epsilon^{*}r_{\caly}^{*}l_{K}^{*}\zeta - l_{K}^{*} \zeta , \
\zeta ) \in T^{*}_{(g, x, z)}(G \times X\times X)$
is  conormal  to $\Omega$.
\end{cor}

\begin{lem}
\label{lem:1.6}
Suppose that $\poidd{G}{P}$ is a Poisson groupoid. Then
for any $f\in C^{\infty}(P)$ and $u\in P$, $X_{\alpha^{*}f}(u)-X_{\beta^{*}f}(u)$
is tangent to $P$ and equals  $X_{f}(u)$.
\end{lem}
\pf As a covector  of $G $ at $u$, $\alpha^{*}df-\beta^{*}df$
is clearly conormal to the unit space $P$. Since $P$ is a
coisotropic submanifold, it follows that $X_{\alpha^{*}f}(u)-X_{\beta^{*}f}(u)$ is tangent to $P$. Hence
$$X_{\alpha^{*}f}(u)-X_{\beta^{*}f}(u)=\alpha_{*} (X_{\alpha^{*}f}(u)-X_{\beta^{*}f}(u))=X_{f}(u) .$$
\qed

\begin{pro}
\label{pro:1.7}
   Suppose that $\poidd{G}{P}$ is a Poisson groupoid acting on
a Poisson manifold $X$ with  moment map $J: X \lon P$.
Suppose that the action is a Poisson action.
Then,   for any $f\in C^{\infty}(P)$,
$$X_{J^{*}f}(x)=(r_{x})_{*} X_{\alpha^{*}f}(u),  \ \forall x\in X$$
where $u=J(x)$.
\end{pro}
\pf Take any $g\in G$ such that $\beta (g)=J(x)$.
Then $(g, x, z)$ with $z=gx$ is in $\Omega $.
  For any $\zeta \in T^{*}_{z}X$ and $f\in C^{\infty}(P)$,
as  covectors at $(g, x, z)$, both
 $(-r_{\caly}^{*}\zeta ,
J^{*}\epsilon^{*}r_{\caly}^{*}l_{K}^{*}\zeta  -l_{K}^{*}\zeta , \zeta )$
and  $(-\beta^{*}df  , J^{*}df  , 0)$
are conormal  to $\Omega$.
Therefore,
$$\pi_{G }(g) (-r_{\caly}^{*}\zeta , -\beta^{*}df )
+\pi (x) (J^{*}\epsilon^{*}r_{\caly}^{*}l_{K}^{*}\zeta
-l_{K}^{*}\zeta ,\  J^{*}df )=0.$$
Now
\be
  \pi_{G }(g) (-r_{\caly}^{*}\zeta , -\beta^{*}df )
&=&
<-X_{\beta^{*}f}(g) , r_{\caly}^{*} \zeta >\\
&=&<-r_{\caly_{*}}X_{\beta^{*}f}(g) ,\zeta>\\
&=&<-r_{\caly_{*}}l_{K_{*}}X_{\beta^{*}f}(u), \zeta >,
\ee
where the last step uses the fact: $X_{\beta^{*}f}(g)=l_{K_{*}}X_{\beta^{*}f}(u)$.

Also,
\be
\pi (x) (J^{*}\epsilon^{*}r_{\caly}^{*}l_{K}^{*}\zeta , J^{*}df )&=&
(J_{*}\pi (x) ) ( \epsilon^{*}r_{\caly}^{*}l_{K}^{*}\zeta , df )\\
&=&\pi_{P}(u)( \epsilon^{*}r_{\caly}^{*}l_{K}^{*}\zeta , df )\\
&=&-X_{f}  (u) (\epsilon^{*}r_{\caly}^{*}l_{K}^{*}\zeta )\\
&=&<-l_{K_{*}}r_{\caly_{*}}\epsilon_{*}X_{f}(u), \zeta > (\mbox{Using Lemma
\ref{lem:1.6})} \\
&=&<-l_{K_{*}}r_{\caly_{*}} (X_{\alpha^{*}f}(u)-
X_{\beta^{*}f}(u)), \zeta >, \\
\ee

and
$$\pi (x) (-l_{K}^{*}\zeta , J^{*}df )=<X_{J^{*}f }(x),  l_{K}^{*}\zeta >
=<l_{K_{*}}X_{J^{*}f}(x), \zeta >.$$
Therefore, it follows that
$<-l_{K_{*}}r_{\caly_{*}}X_{\alpha^{*}f}(u), \zeta >+<l_{K_{*}}X_{J^{*}f}(x),
 \zeta >=0$, which implies immediately that
$$X_{J^{*}f}(x)=r_{\caly_{*}}X_{\alpha^{*}f}(u)=(r_{x})_{*}X_{\alpha^{*}f}(u).$$\qed

\begin{pro}
\label{pro:multi}
   Suppose that $\poidd{G}{P}$ is a Poisson groupoid acting on
a Poisson manifold $X$ with  moment map $J: X \lon P$.
Suppose that the action is a Poisson action. Then,
\begin{equation}
 \pi (gx)=(l_{K })_{*} \pi (x) + (r_{\caly })_{*} \pi_{G } (g)
- (r_{\caly })_{*}(l_{K })_{*}\pi_{G}(u),
 \end{equation}
where $u=\beta (g)=J(x)$, $K $ is any bisection of
$G$ through $g$, and $\caly$ any local section through
the point $x$.
\end{pro}

\begin{lem}
For any $u\in P$ and $\eta  \in T^{*}_{u}G$,
\begin{equation}
X_{\alpha^{*}\epsilon^{*} \eta }(u)+\epsilon_{*}\alpha_{*}X_{\eta}(u)
-\epsilon_{*}X_{\epsilon^{*}\eta }(u)=X_{\eta }(u) .
\end{equation}
\end{lem}
\pf  It is clear that $\alpha^{*}\epsilon^{*}\eta -\eta$, as
a covector of $G$ at $u$, is conormal to $P$.
Hence, $X_{\alpha^{*}\epsilon^{*} \eta }(u)-X_{\eta }(u)$
is tangent to $P$. It follows that
\be
X_{\alpha^{*}\epsilon^{*} \eta }(u)-X_{\eta }(u)&=&
\epsilon_{*}\alpha_{*}(X_{\alpha^{*}\epsilon^{*} \eta }(u)-X_{\eta }(u))\\
&=&\epsilon_{*}X_{\epsilon^{*}\eta }(u)-\epsilon_{*}\alpha_{*}X_{\eta}(u).
\ee
This completes the proof of the lemma. \qed

\begin{lem}
\label{lem:1.9}
$$\pi (x) (J^{*}\epsilon^{*}r_{\caly}^{*}l_{K}^{*}\zeta  -l_{K}^{*}\zeta ,
J^{*}\epsilon^{*}r_{\caly}^{*}l_{K}^{*}\eta  -l_{K}^{*}\eta )=
(l_{K_{*}}\pi (x) )(\zeta , \eta )-(l_{K_{*}}r_{\caly_{*}}\pi_{G }(u))
(\zeta , \eta ).$$
\end{lem}
\pf
\be
\pi (x) (J^{*}\epsilon^{*}r_{\caly}^{*}l_{K}^{*}\zeta , J^{*}\epsilon^{*}r_{\caly}^{*}l_{K}^{*}\eta ) &=&
(J_{*}\pi (x))(\epsilon^{*}r_{\caly}^{*}l_{K}^{*}\zeta , \epsilon^{*}r_{\caly}^{*}l_{K}^{*}\eta )\\
&=&\pi_{P}(u)(\epsilon^{*}r_{\caly}^{*}l_{K}^{*}\zeta , \epsilon^{*}r_{\caly}^{*}l_{K}^{*}\eta )\\
&=&<X_{\epsilon^{*}r_{\caly}^{*}l_{K}^{*}\zeta} (u),  \epsilon^{*}r_{\caly}^{*}l_{K}^{*}\eta >\\
&=&<\epsilon_{*} X_{\epsilon^{*}r_{\caly}^{*}l_{K}^{*}\zeta} (u) ,
r_{\caly}^{*}l_{K}^{*}\eta >.
\ee

\be
\pi (x) (J^{*}\epsilon^{*}r_{\caly}^{*}l_{K}^{*}\zeta ,\  l_{K}^{*}\eta )&=&
<X_{J^{*}\epsilon^{*}r_{\caly}^{*}l_{K}^{*}\zeta }(x), l_{K}^{*}\eta >
\mbox{(  using Proposition \ref{pro:1.7} )} \\
&=&< r_{\caly_{*}}X_{\alpha^{*}\epsilon^{*}r_{\caly}^{*}l_{K}^{*}\zeta }(u),
l_{K}^{*}\eta >\\
&=& <X_{\alpha^{*}\epsilon^{*}r_{\caly}^{*}l_{K}^{*}\zeta }(u) , r_{\caly}^{*}l_{K}^{*}\eta >.
\ee
Using the above relations,

\be
\pi (x) ( l_{K}^{*}\zeta , J^{*}\epsilon^{*}r_{\caly}^{*}l_{K}^{*}\eta )&=&
-< X_{\alpha^{*}\epsilon^{*}r_{\caly}^{*}l_{K}^{*}\eta }(u) , r_{\caly}^{*}l_{K}^{*}\zeta >\\
&=& <X_{r_{\caly}^{*}l_{K}^{*}\zeta }(u),   \alpha^{*}\epsilon^{*}r_{\caly}^{*}l_{K}^{*}\eta >\\
&=&<\epsilon_{*}\alpha_{*}X_{r_{\caly}^{*}l_{K}^{*}\zeta }(u), r_{\caly}^{*}l_{K}^{*}\eta >.
\ee

\be
&&\pi (x) (J^{*}\epsilon^{*}r_{\caly}^{*}l_{K}^{*}\zeta  -l_{K}^{*}\zeta ,
J^{*}\epsilon^{*}r_{\caly}^{*}l_{K}^{*}\eta  -l_{K}^{*}\eta  )\\
&=&
\pi (x) (J^{*}\epsilon^{*}r_{\caly}^{*}l_{K}^{*}\zeta , J^{*}\epsilon^{*}r_{\caly}^{*}l_{K}^{*}\eta ) -
\pi (x) (J^{*}\epsilon^{*}r_{\caly}^{*}l_{K}^{*}\zeta , l_{K}^{*}\eta )\\
& & -
\pi (x) ( l_{K}^{*}\zeta , J^{*}\epsilon^{*}r_{\caly}^{*}l_{K}^{*}\eta )
+\pi (x) (l_{K}^{*}\zeta , l_{K}^{*}\eta )\\
&=&<\epsilon_{*}X_{\epsilon^{*}r_{\caly}^{*}l_{K}^{*}\zeta }(u)
-X_{\alpha^{*}\epsilon^{*}r_{\caly}^{*}l_{K}^{*}\zeta }(u)
-\epsilon_{*}\alpha_{*}X_{r_{\caly}^{*}l_{K}^{*}\zeta }(u) , r_{\caly}^{*}l_{K}^{*}\eta >  +\pi (x) (l_{K}^{*}\zeta , l_{K}^{*} \eta )\\
&=&<-X_{r_{\caly}^{*}l_{K}^{*}\zeta}(u), r_{\caly}^{*}l_{K}^{*}\eta >
+\pi (x) (l_{K}^{*}\zeta , l_{K}^{*} \eta )\\
&=&-\pi_{G }(u) (r_{\caly}^{*}l_{K}^{*}\zeta , r_{\caly}^{*}l_{K}^{*}\eta )
+\pi (x) (l_{K}^{*}\zeta , l_{K}^{*} \eta )\\
&=&-(l_{K_{*}} r_{\caly_{*}}\pi_{G }(u) (\zeta , \eta ) +(l_{K_{*}}\pi (x) )(\zeta , \eta )
\ee
\qed
{\bf Proof of Proposition \ref{pro:multi}: } For any $\zeta , \eta \in T^{*}_{z}X$,
both  $(-r_{\caly}^{*}\zeta , J^{*}\epsilon^{*}r_{\caly}^{*}l_{K}^{*}\zeta - l_{K}^{*} \zeta ,
\zeta )$  and $(-r_{\caly}^{*}\eta , J^{*}\epsilon^{*}r_{\caly}^{*}l_{K}^{*}\eta - l_{K}^{*} \eta ,
\eta )$  are conormal  to $\Omega $. Therefore,

$$\pi_{G }(g)(-r_{\caly}^{*}\zeta , -r_{\caly}^{*}\eta ) +
\pi (x) (J^{*}\epsilon^{*}r_{\caly}^{*}l_{K}^{*}\zeta  -l_{K}^{*}\zeta ,
J^{*}\epsilon^{*}r_{\caly}^{*}l_{K}^{*}\eta  -l_{K}^{*}\eta )
-\pi (z) (\zeta , \eta )=0.  $$
The conclusion follows immediately by Lemma \ref{lem:1.9}. \qed
{\bf Proof of Theorem \ref{thm:main1}: }  One
direction has been proved by the above series of propositions.
It remains to prove the other direction.

First, we note that the first condition implies that $J: X\lon P$ is a Poisson map.
This can be seen as follows.

 As a map defined on  $\beta^{-1}(u)$,
we have  $(J\smalcirc r_{x})(r)=J(rx)=\alpha (r)$, $\forall r\in \beta^{-1}(u)$.
Then, it follows that
$$J_{*}X_{J^{*}f}(x)=J_{*} (r_{x})_{*}X_{\alpha^{*}f}(u)
=\alpha_{*}X_{\alpha^{*}f}(u)=X_{f}(u).$$

Given  any point $(g, x , z)\in \Omega \subset G \times  X\times \overline{X}$.
The   conormal  space of $\Omega$ at this point   is spanned by two types of vectors:
$(-\beta^{*}df  , J^{*}df  , 0)$ for any $f\in C^{\infty}(P)$ and
$(-r_{\caly}^{*}\zeta , J^{*}\epsilon^{*}r_{\caly}^{*}
l_{K}^{*}\zeta - l_{K}^{*} \zeta , \zeta )$  for any $\zeta \in T_{z}^{*}X$.
Using the  same arguments as in the
proof of  Propositions \ref{pro:1.2}, 
\ref{pro:1.7} and \ref{pro:multi}, it can be
easily checked
 that the evaluation of the Poisson tensor on all these
vectors vanish. This concludes the proof. \qed
{\bf Remark.}  In the proof above,  we
only used  the fact that  Equation (\ref{eq:multi})
holds for one, instead of all, such bisections. 
Consequently, under the rest of  the assumptions of Theorem
\ref{thm:main1},
if Equation (\ref{eq:multi}) holds for any one  bisection,
it holds for all.

\section{Poisson homogeneous spaces}

The notion of homogeneous space for a groupoid action is more subtle
than for groups. (This point has already been made in
\cite{br-da-ha:topological}, cited in \cite{Mackenzie:book}). One
natural candidate for such a space is $G$ acting on itself by left
translations, but this action is not transitive in the usual sense,
since $\beta(gx)=\beta(x)$, so that the action is transitive only on
each $\beta$-fibre.

Instead, we define homogeneous $G$-spaces to be those which are
isomorphic to $G/H$ for some wide (i.e.
containing all the identities) subgroupoid
$H$ of
$G$.  That is, we define $G/H$ by the equivalence relation $g\sim h
\iff \exists n \in H $ such that $gn=h$, with the moment map $J([g])=\alpha(g)$
and the action $g\cdot [h] = [gh]$.  The following is an intrinsic
characterization of such spaces.

\begin{defi}
\label{defi-homogeneous}
A $G$-space $X$ over $P$ is {\bf homogeneous} if there is a section $\sigma$
of the moment map $J$ which is {\bf saturating} for the action in the sense that
$G
\cdot\sigma(P) = X$.  The {\bf isotropy subgroupoid} of the section $\sigma$
consists of those
$g\in G$ for which $g\cdot \sigma(P) \subset \sigma(P)$.
\end{defi}

\begin{pro}
\label{pro-homogeneous}
A $G$-space is homogeneous if and only if it is isomorphic to $G/H$ for
some wide subgroupoid $H\subset G$.
\end{pro}
\pf It is easy to see that a $G$-space of the form $G/H$ is homogeneous,
since the image of the identity section of $G$ is a saturating section of
$G/H$.   On the other hand, given a homogeneous space $X$ with saturating
section
$\sigma$, we define a map $\theta:G\lon X$  by
$\theta(g)=g\cdot\sigma(\beta (g))$.  $\theta$ is surjective: for every $x\in X$
we have $x=g\cdot \sigma(p)$ for some $g\in G$ and $p\in P$; then
$\beta(g)=J(\sigma(p))=p$, so $x=g\cdot\sigma (\beta (g))=\theta(g).$  On the
other hand, if
$\theta(g)=\theta(h)$, i.e. $g\cdot\sigma(\beta (g))=h\cdot\sigma(\beta (h))$,
we have $g^{-1}h\cdot \sigma(\beta (h)) = \sigma(\beta (g))$.  Since
$g^{-1}h$ can act on only one element of $\sigma(P)$, it follows that $g^{-1}h$
belongs to the isotropy groupoid $H$ of the section $\sigma$.  Also, for $n\in
H$, $\theta(gn)=gn\cdot \sigma(\beta(gn))=g\cdot\sigma(p)$ for some $p$.  But
$\beta(g)=J(\sigma(p))$, so  $$\theta(gn)=g\cdot(\sigma(\beta(g))=\theta(g).$$ It
follows that
$\theta$ induces a bijection (which can be checked to be
$G$-equivariant) from $G/H$ to $X$.
\qed

Let $\poiddd{G}{P}{}$ be a  Poisson groupoid, and  $H$ a connected
closed  subgroupoid.
 Suppose  that $X=G /H$ is  a Poisson manifold. Write $p$  as
the  natural projection: $G \lon X$.
 Let $\{\cdot , \cdot\}_{1}$ be the  difference
bracket from $C^{\infty}(X) \otimes C^{\infty}(X)$
to $C^{\infty}(G  )$ as defined  by Equation (\ref{eq:b1}):
 $$\{\F,\G\}_1  = p^*\{\F,\ \G\} - \{p^* \F,\ p^* \G\}_{G }, 
 \ \ \ \forall \F,\ \G\in   C^{\infty}(X).$$

$X$ is a homogeneous space of the groupoid $G $. Write $J$ as
its moment map: $X\lon P$. Then, $J\smalcirc p=\alpha$.
The main theorem of this section is

\begin{thm}
\label{thm:main2}
$X$ is a Poisson homogeneous space iff
 \begin{enumerate}
\item  For any bisection $K$ of $G$,
$$\{l_{K}^{*}\F,\  l_{K}^{*}\G\}_{1}=l_{K}^{*}\{\F, \ \G\}_{1}, \ \ \ \forall
\F, \G\in C^{\infty}(X), $$
i.e.,  $\{\cdot , \cdot \}_1$ is left invariant; and
\item for any $f\in C^{\infty}(P)$ and $\G\in C^{\infty}(X)$,
   $$\{J^*f, \G \}_1=0. $$
\end{enumerate}
\end{thm}

We split its  proof 
 into two propositions.

\begin{pro}
The following statements are equivalent:
\begin{enumerate}
\item For any $f\in C^{\infty}(P)$,  $X_{J^{*}f}(x)=(r_{x})_{*}X_{\alpha^{*}f}(u)$, $\forall x\in X$ with $u=J(x)$.
\item   For any $f\in C^{\infty}(P)$ and $\G\in C^{\infty}(X)$
   $\{J^*f, \  \G\}_1=0$.
\end{enumerate}
\end{pro}
\pf For any $g\in G$, let $x=p (g)=[g H] \in X$. Then, $\forall f\in
C^{\infty} (P)$,
\be
\{p^{*}J^{*}f, \ \ p^{*}\G\}_{G }(g)&=& \{\alpha^{*}f,\ \  p^{*}\G\}_{G}(g)\\
 &=&X_{\alpha^{*}f}(g)( p^{*}\G)\\
&=&[p_{*} X_{\alpha^{*}f}(g)]\G\\
&=&[p_{*} r_{g_{*}} X_{\alpha^{*}f}(u)]\G
\ee
Now, it is clear that
$$(p\smalcirc  r_{g})(\gamma )=p (\gamma g)=[\gamma gH]=r_{x} (\gamma ), \ \ \ \
\forall \gamma \in \beta^{-1} (u). $$
Hence, $p\smalcirc  r_{g}=r_{x}$ on $\beta^{-1}(u)$.
Therefore, $$\{p^{*}J^{*}f,\ \  p^{*}\G\}_{G }(g)=
[(r_{x})_{*}X_{\alpha^{*}f}(u)]\G.  $$
On the other hand,

\be
p^{*}\{J^{*}f ,\ \  \G\}(g)&=&\{J^{*}f, \ \G\}(pg)\\
&=&\{J^{*}f, \ \G\}(x)\\
&=&X_{J^{*}f}(x)\G.
\ee

Thus, it follows immediately that
$\{J^*f,\  \G\}_1=0$, $\forall \G\in C^{\infty}(X)$ is equivalent
to the equation:
  $X_{J^{*}f}(x)=(r_{x})_{*}X_{\alpha^{*}f}(u)$, $\forall x\in X$.  \qed

\begin{pro}
The following  statements are equivalent:
\begin{enumerate}
\item   For any bisection $K$ of $G$,
$$\{l_{K}^{*}\F,\  l_{K}^{*}\G\}_{1}=l_{K}^{*}\{\F, \ \G\}_{1}, \ \ \ \forall
\F, \G\in C^{\infty}(X), $$
i.e.,  $\{\cdot , \cdot \}_1$ is left invariant;
\item   For any compatible $(g, x)\in G *X$,
\begin{equation}
 \pi (gx)=(l_{K })_{*} \pi (x) + (r_{\caly })_{*} \pi_{G } (g)
- (r_{\caly })_{*}(l_{K })_{*}\pi_{G}(u),
 \end{equation}
where $u=\beta (g)=J(x)$, $K $ is any local bisection of
$G$ through $g$, and $\caly$ any local section through
the point $x$.
\end{enumerate}
\end{pro}
\pf  Let $\gamma $ be any point in $G$,
$p (\gamma )=[\gamma H ]=x\in
X$,  and $J(x)=\alpha (\gamma )=u\in P$.
Also, let $g=K( u)=K\cap \beta^{-1}(u)\in G $.
Then, $l_{K}\gamma =g \gamma $ and $l_{K }x=g\cdot x$. Thus,

\be
p^{*}\{l_{K}^{*}\F, \  l_{K}^{*}\G\}(\gamma )&=&\{l_{K}^{*}\F,  \ l_{K}^{*}\G\}(x)\\
&=&\pi (x) (l_{K}^{*}d\F , \  l_{K}^{*}d\G )\\
&=&l_{K_{*}}\pi (x)(d\F,  \ d\G),
\ee
and

\be
&&\{p^{*}l_{K}^{*}\F,\  p^{*}l_{K}^{*}\G\}_{G }(\gamma )\\
&=&\pi_{G} (\gamma )(p^{*}l_{K}^{*}d\F, \ p^{*}l_{K}^{*}d\G )\\
&=&l_{K_{*}}p_{*} \pi_{G} (\gamma ) (d\F, \ d\G).
\ee

That is,

$$\{l_{K}^{*}\F, \ l_{K}^{*}\G\}_{1}(\gamma )=-l_{K_{*}}\pi (x)(d\F, d\G)
+l_{K_{*}}p_{*} \pi_{G} (\gamma ) (d\F, d\G) .$$

On the other hand,

\be
l_{K}^{*}\{\F ,\G\}_{1}(\gamma )&=&\{\F, \G\}_{1}(g\gamma ) \\
&=&-(p^{*}\{\F, \G\})(g\gamma )+\{p^{*}\F, p^{*}\G\}_{G }(g \gamma )\\
&=&-\{\F, \G\}(gx) +\{p^{*}\F ,p^{*}\G\}_{G }(g\gamma )\\
&=&-\pi (gx)(d\F ,d\G )+p_{*}\pi_{G }(g\gamma )(d\F, d\G).
\ee

Therefore, $\{\cdot , \cdot \}_{1} $ is left invariant iff
$$l_{K_{*}}\pi (x)-l_{K_{*}}p_{*}\pi_{G} (\gamma )=
\pi (gx) -p_{*}\pi_{G }(g\gamma ), $$ or

\begin{equation}
\label{eq:multi-eq}
p_{*}\pi_{G }(g\gamma )-l_{K_{*}}p_{*}\pi_{G} (\gamma )=
\pi (gx) -l_{K_{*}}\pi (x).
\end{equation}

Since $G $ is a  Poisson groupoid, according to Theorem 2.4 in \cite{Xu},
we have
$$ \pi_{G} (g \gamma )=l_{K_{*} } \pi_{G}  (\gamma ) + r_{\calr_{*} }
\pi_{G } (g)- l_{K_{*} } r_{\calr_{*} }  \pi_{G}(u), $$
where $\calr$ is any local  bisection of $G $
 throuhg the point $\gamma $.
Hence, it follows that
\be
 p_{*}\pi_{G }(g\gamma )-l_{K_{*}}p_{*}\pi_{G} (\gamma ) &=&
p_{*}r_{\calr_{*}}\pi_{G }(g) -l_{K_{*}}p_{*}r_{\calr_{*}}\pi_{G }(u).
\ee
Here, we have used the identity: $p\smalcirc l_{K}=l_{K}\smalcirc p$,
as both being considered as maps from $G $ to $X$.

Let $\caly =p (\calr )\subset X$. Then $\caly$ is a  local section
of $J$ through the point $x$.
It is simple to see that,
as maps from $G $ to $X$,  $p\smalcirc r_{\calr} =r_{\caly}$.

Hence,
$$p_{*}r_{\calr_{*}}\pi_{G }(g) -l_{K_{*}}p_{*}r_{\calr_{*}}\pi_{G }(u)=
r_{\caly_{*}}\pi_{G }(g) -l_{K_{*}}r_{\caly_{*}}\pi_{G}(u).$$

This shows that Equation (\ref{eq:multi-eq}) is
equivalent to  Equation (\ref{eq:multi}) with $\caly=p  (\calr )$.
 The conclusion thus follows from
Theorem \ref{thm:main1} together with
 the remark following its proof. \qed


A Dirac structure $L$ of $A \oplus A^*$ is called \closed if $L \cap A$
is a subalgebroid of $A$ whose left translation
 defines a simple foliation on $G$\footnote{In the case of groups,
this is equivalent to saying that $L\cap A$ can be
integrated to a connected closed subgroup of $G$. However, when
$G$ is a groupoid, that $L\cap A$ can be integrated to a
connected closed subgroupoid 
 seems not sufficient to get a simple
foliation.}.
Then,   $L$ is \closed iff  $\bar L$ is \regular.
Thus,  combining Theorem
\ref{thm:main4} and  Theorem \ref{thm:main2}, we obtain
the following main theorem,  which is a generalization
of Drinfel'd's theorem in the groupoid context.

\begin{thm}
\label{thm:main5}
For  a Poisson groupoid $G$, 
there is an 1-1 correspondence between   Poisson homogeneous spaces
$G / H$ and  \closed Dirac structures $L$
 of its tangent Lie bialgebroid, where $H$ is the $\alpha$-connected closed
subgroupoid of $G$ corresponding to the subalgebroid  $L \cap A$.
\end{thm}

We end this section with some examples.

\begin{ex}
\label{ex:hom1}
Under the same hypothesis  as  in Theorem \ref{thm:main5},
if moreover $L$ is 
the graph of a   hamiltonian
operator $ H\in \gm (\wedge^{2} A ) $, its
corresponding Poisson homogeneous space is 
still $G$, but equipped with a different
Poisson structure  $\pi_G + \Phi^*H $. Here 
  $\Phi^*H \in \gm (\wedge^{2}  TG) $ is  the pull back
of $H$ under the morphism $\Phi : T^{*} G\lon A^*$.
This  Poisson structure in fact defines 
 a  Poisson affinoid structure in terms of Weinstein \cite{we:affinoid}.
\end{ex}

\begin{ex}
\label{ex:hom2}
For a Poisson manifold $(P, \pi )$,
let  $(TP, T^* P,\pi)$ be  its  canonical
Lie bialgebroid.
The  corresponding Poisson groupoid is  the pair groupoid
$G =  P \times {\bar P}$.
It is easy to see that  any  homogeneous $G$-space in this case
 is always of the form
$ P \times P/{\cal F}$,  where ${\cal F}$ is a simple foliation on $P$.
The groupoid $G$-action  is given by:
$$ (x, y) \cdot (y, [z]) = (x, [z]),  \ \  \ \ \ \forall x, y, z \in P.$$
Moreover, this becomes  a Poisson homogeneous $G$-space iff  $P$ is equipped
with the original Poisson structure $\pi$ (the Poisson structure
on $ P/{\cal F}$ may be  arbitrary). In other words, in this case,
 Poisson homogeneous
spaces are in 1-1 correspondence with Poisson  structures on 
 quotient manifolds of $P$.
Thus, Theorem \ref{thm:main5} reduces to Theorem  \ref{thm:main02}.
\end{ex}

\begin{ex}
\label{ex:hom3}
Dually, we may   switch the order  and
 consider the Lie bialgebroid $(T^{*}P, TP )$   for 
 a Poisson manifold  $P$.
Its  corresponding Poisson groupoid $G$, if it exists,  is in fact a  symplectic
groupoid of  $P$ (see Theorem 5.3 in \cite{MackenzieX:1996}).
  It is not difficult to see that
   a  homogeneous space  $X$  becomes a   Poisson homogeneous space iff
the moment map $J: X\lon P$ is a  Poisson map.
Thus we obtain the following 

\begin{cor}
Suppose that $P$ is an integrable  Poisson manifold with symplectic
groupoid $G$. There is a one-one
correspondence between \regular Dirac structures in
the double  $E=T^{*}P \oplus TP$ and
 homogeneous  $G$-spaces $X$ equipped with   a compatible
Poisson structure in the sense that the moment
map $J: X\lon P$ is a  Poisson map.
\end{cor}
\end{ex}

It is worth noting that when $P$ is symplectic,  for
a given null  Dirac structure  on $E= TP \oplus T^{*}P$,
the  corresponding pair of Poisson   homogeneous spaces
in Example \ref{ex:hom2} and \ref{ex:hom3} correspond to a
Poisson dual pair. In other words, for a symplectic
manifold $P$,   a null  Dirac structure on $E= TP \oplus T^{*}P$,
under a certain regularity condition, corresponds to
a  Poisson dual pair.
It would be interesting to explore what happens for a
general Dirac structure, and also even   more  general
situation when $P$ is degenerate.

 The last example is the following

\begin{ex}
As in Example \ref{ex:hom3}, let $P$ be an integrable Poisson manifold
with symplectic groupoid $G$, and
$E=T^{*} P \oplus TP$.

Assume that the Dirac structure arises from a hamiltonian operator,
which is, in this case, a two form $\theta $ on $P$ satisfying
the equation:

$$d\theta +\half [\theta , \theta ]=0.$$
Its corresponding Poisson homogeneous space, as described in
Example \ref{ex:hom1}, is $G$ equipped with the ``affine"
Poisson structure  $\pi_{G} + \Phi^{*}\theta$, where
$\Phi:  T^* G \lon TP$ is the $\beta$-map of the cotangent
symplectic groupoid as defined by Equation (\ref{eq:Phi}).
In general, it is not clear whether this is still nondegenerate.

However, in the extreme case that $P$ is a zero Poisson
structure, we will see that it is still symplectic.

To see this, we
note that $\Phi$ 
fits into the following commutative  diagram:

\begin{equation}                         
\label{eq:diagram}
\matrix{&\Phi&\cr
        T^* G &\vlra&TP\cr
        &&\cr \pi^{\#}_G
\Bigg\downarrow &&\Bigg\downarrow id\cr
        &&\cr
        TG&\vlra&TP\cr
        &\beta_* &\cr}.
\end{equation}
This implies that
$$(\Phi^*\theta)^{\#}=-\pi_{G}^{\#}\smalcirc (\beta^{*}\theta )^{b}
\smalcirc \pi_{G}^{\#}. $$
Therefore, 
$$(\Phi^*\theta)^{\#}\smalcirc \omega^{b}=-\pi_{G}^{\#}\smalcirc (\beta^{*}\theta )^{b},$$
where $\omega$ denotes the symplectic structure on $G$.

On the other hand, it follows from the
fact that $Im(\alpha^*\theta)^{b} \subset T^{\alpha}G^{\perp}$ and
$T^{\alpha}G^{\perp} \subset ker\Phi$ that 
$$(\Phi^*\theta)^{\#}\smalcirc   (\alpha^*\theta)^{b} =0. $$

Thus,
$$(\pi_G^{\#} + (\Phi^{*}\theta )^{\#}) \smalcirc ( \omega^{b} + (\alpha^{*}\theta)^{b}) = 
id + \pi_G^{\#} \smalcirc ((\alpha^{*}\theta)^{b} - (\beta^{*}\theta)^{b} ).$$

That is, in case that $\beta^*\theta = \alpha^*\theta $, the affine
Poisson structure  $\pi_G + \Phi^*\theta$ is non-degenerate and the corresponding
symplectic form is $\omega + \alpha^*\theta$.

Thus, when  $P$ is a  zero Poisson  manifold, 
its symplectic groupoid  is
$T^{*}P$ equipped  with the 
canonical cotangent symplectic structure.
In this case, $\alpha = \beta$ and is just
 the natural projection from $T^*P$ to $P$. A  hamiltonian
operator corresponds to any closed two form $\theta$ on $P$.
The homogeneous space corresponding to its induced Dirac
structure is 
again $T^{*}P$, with the non-degenerate Poisson structure coming from 
the sum of the canonical 2-form and the pullback of $\theta $ by the 
projection $T^{*}P\lon P$.
\end{ex}

     \end{document}